\documentclass[10pt,a4paper]{article}

\usepackage{geometry}
\usepackage[utf8]{inputenc}
\usepackage{amsthm, amssymb, natbib, graphicx}
\usepackage{listings}
\usepackage[ruled,vlined]{algorithm2e}
\usepackage{amsmath,tikz}
\usepackage{array, makecell}
\usepackage{sidecap, caption}
\usepackage{xcolor, color} 
\usepackage{comment} 
\usepackage{algorithm2e}
\usepackage{caption}
\usepackage{subcaption}
\usepackage[most]{tcolorbox}
\setcitestyle{numbers, comma,square}
\usepackage{bbm}
\usepackage{tabularray}
\usepackage{diagbox}
\UseTblrLibrary{booktabs}

\newtheorem{theorem}{Theorem}[section]
\newtheorem{corollary}{Corollary}[section]
\newtheorem{lemma}[theorem]{Lemma}

\newtheorem{remark}{Remark}[section]
\newtheorem{assumption}{Assumption}
\newtheorem{definition}{Definition}[section]
\DeclareMathOperator*{\argmin}{argmin}

\newcommand{\R}{\mathbb{R}}

\usepackage[normalem]{ulem}

\newcommand{\bb}[1]{\textcolor{blue}{#1}}

\newcommand{\ignore}[1]{}
\def\T{{\cal T}}
\def\C{{\cal C}}

\def\M{{\cal{M}}}
\def\N{{\cal{N}}}

\title{Global Convergence of High-Order Regularization Methods with Sums-of-Squares Taylor Models}
\author{\thanks{wenqi.zhu@maths.ox.ac.uk, University of Oxford, UK} Wenqi Zhu,  \thanks{cartis@maths.ox.ac.uk, University of Oxford, UK} Coralia Cartis}

\author{Coralia Cartis\thanks{The order of the authors is alphabetical; the second author (Wenqi Zhu) is the primary contributor.  This work was supported by the Hong Kong Innovation and Technology Commission (InnoHK Project CIMDA).} \textsuperscript{\normalfont ,}\thanks{Mathematical Institute, Woodstock Road, University of Oxford, Oxford, UK, OX2 6GG.  {\tt cartis@maths.ox.ac.uk} } \quad and \quad
Wenqi Zhu\footnotemark[1] \textsuperscript{\normalfont,}\thanks{Mathematical Institute, Woodstock Road, University of Oxford, Oxford, UK, OX2 6GG.
{\tt wenqi.zhu@maths.ox.ac.uk}}}

\begin{document}

\maketitle

\begin{abstract}
High-order tensor methods that employ  Taylor-based local models  (of degree  $p\ge 3$) within adaptive regularization frameworks have been recently proposed for both convex and nonconvex optimization problems.  They have been shown to have superior, and even optimal, worst-case global convergence rates and local rates compared to Newton's method. 
Finding rigorous and efficient techniques for minimizing the Taylor polynomial sub-problems remains a challenging aspect 
for these algorithms.
Ahmadi et al \cite{ahmadi2023higher} recently introduced a tensor method based on sum-of-squares (SoS) reformulations, so that each Taylor polynomial sub-problem in their approach can be tractably minimized using semidefinite programming (SDP) \cite{ahmadi2023higher}; however, the global convergence and complexity of their method have not been addressed for general nonconvex problems.
This paper introduces an algorithmic framework that combines the Sum of Squares (SoS) Taylor model with adaptive regularization techniques for nonconvex smooth optimization problems. Each iteration minimizes an SoS Taylor model, offering a polynomial cost per iteration. For general nonconvex functions, the worst-case evaluation complexity bound is $\mathcal{O}(\epsilon^{-2})$, while for strongly convex functions, an improved evaluation complexity bound of $\mathcal{O}(\epsilon^{-\frac{1}{p}})$ is established. To the best of our knowledge, this is the first global rate analysis for an adaptive regularization algorithm with a tractable high-order sub-problem in nonconvex smooth optimization, opening the way for further improvements.
\end{abstract}

\section{Introduction}
In this paper, we consider the unconstrained nonconvex optimization problem, 
\begin{equation}
  \min_{x\in \R^n} f(x)  
  \label{min f}
\end{equation}
where  $f : \R^n \rightarrow \R $ and $f  \in \C^{p, 1}(\R^n) $ which means that $f$ is $p$-times continuously differentiable ($p \ge 1$), bounded below and the $p$th derivative of $f$ is globally Lipschitz continuous. Recent research, including papers by Birgin et al. \cite{birgin2017worst}, Cartis et al. \cite{cartis2022evaluation, cartis2020sharp, cartis2020concise}, and Carmon et al. \cite{carmon2020lower, carmon2021lower}, has shown that certain optimization algorithms exhibit superior worst-case global complexity bounds even for potentially nonconvex objective functions. 
In these methods \cite{cartis2022evaluation, cartis2020sharp, cartis2020concise}, a polynomial local model $m_p(x, s)$ is typically constructed to approximate $f(x+s)$ at the iterate $x=x_k$. Then, $x_k$ is iteratively updated by $
s_k \approx \argmin_{s\in \R^n}{m_p(x_k, s)}$, and $x_{k+1} := x_k+s_k$, whenever sufficient objective decrease is obtained, until an approximate local minimizer of $f$ is generated. 
The minimizer approximately satisfies the first-order optimality condition\footnote{Unless otherwise stated, $\|\cdot\|$ denotes the Euclidean norm in this paper.}, represented by $\|\nabla_x f(x_k)\|\le \epsilon$, where $\epsilon$ represents the prescribed first order tolerance. The construction of the model $m_p$ relies on the $p$th-order Taylor expansion of $f(x_k+s)$ at $x_k$
\begin{equation}
\tag{$p$th Order Taylor Model}
   T_{p}(x_k, s) := f(x_k) + \sum_{j=1}^p \frac{1}{j!} \nabla_x^j f(x_k)[s]^j,
\end{equation}
where $\nabla_x^j f(x_k) \in \R^{n^j}$ is a $j$th-order tensor and $\nabla_x^j f(x_k)[s]^j$ is the $j$th-order derivative of $f$ at $x_k$ along $s \in \R^n$. 
To ensure a bounded below local model $m_p$, an adaptive regularization term\footnote{Unless otherwise stated, the regularization term uses the Euclidean norm.} $\sigma_k > 0$ is added to $T_p$. The model  $m_p$ is 
\begin{equation}
    m_p(x_k, s) = T_p(x_k, s)+ \frac{\sigma_k}{p+1}\|s\|^{p+1}.
    \label{subprob}
\end{equation}
The sub-problem \eqref{subprob} corresponds to the local model in the  $p$th order adaptive regularization (i.e., AR$p$) algorithmic framework~\cite{birgin2017worst, cartis2020sharp, cartis2020concise}. In AR$p$ algorithmic framework, the parameter $\sigma_k$ is adjusted adaptively to ensure progress towards optimality over the iterations. Under Lipschitz continuity assumptions on $\nabla_x^p f(x)$, AR$p$ Algorithm requires no more than $\mathcal{O}(\epsilon^{-\frac{p+1}{p}})$ evaluations of $f$ and its derivatives to compute a first-order local minimizer to accuracy $\epsilon$. 
This result demonstrates that as we increase the order $p$, the evaluation complexity bound improves. This superior, and even optimal \cite{carmon2020lower, cartis2022evaluation}, worst-case complexity bound motivates us to develop an efficient algorithmic implementation for the AR$p$ method. 

The case of $p=1$, \eqref{subprob} corresponds to the steepest descent model, while $p=2$ gives rise to the Adaptive Cubic Regularization model (i.e., ARC or AR$2$ model)\cite{cartis2007adaptive, cartis2011adaptive}. Efficient algorithms have been developed for minimizing the AR$p$ sub-problem for the cases of $p = 2$ \cite{cartis2007adaptive, cartis2011adaptive}. However, for the higher-order AR$p$ method ($p \ge 3$), the efficient minimization of the AR$p$ sub-problem remains an open question.
It has been proven in~\cite{burachik2021steklov, luo2010semidefinite} that finding the global minimum of the AR$3$ sub-problem is NP-hard. Moreover, finding the global minimum of general even-degree polynomials with degree at least $4$ is also NP-hard \cite{ahmadi2023higher, ahmadi2022complexity, murty1987some}.
As a result, researchers have focused on iterative algorithms for locating local minima of the AR$p$ sub-problem for $p \ge 3$. According to \cite{carmon2020lower, cartis2020concise}, finding a local minimum of the AR$p$ sub-problem is sufficient to ensure good complexity for AR$3$.

Nesterov et al. have contributed extensively to the field of convex higher-order methods~\cite{Nesterov2021implementable, Nesterov2020inexact, Nesterov2021inexact, Nesterov2021superfast, Nesterov2022quartic, Nesterov2006cubic}. These methods leverage various convex optimization tools, including Bregman gradient methods, high-order proximal-point operators, and iterative minimization of convex quadratic models with quartic regularization.
While Nesterov's work primarily focuses on convex objectives, recent developments by Cartis and Zhu have introduced the Quadratic Quartic Regularization (QQR) method for the minimization of generally nonconvex AR$3$ sub-problems \cite{cartis2023second, zhu2022quartic}. QQR approximates the third-order tensor term of the AR$3$ model with a linear combination of quadratic and quartic terms. Following the QQR method, the Cubic-Quartic Regularization (CQR) algorithmic framework \cite{zhu2023cubic} was also introduced for minimizing nonconvex cubic multivariate polynomials with quartic regularization. CQR methods leverage the cubic term to approximate local tensor information, while the adaptive quartic term provides model regularization and controls progress.

Recently, Ahmadi et al. \cite{ahmadi2023higher} introduced an iterative algorithm based on sum-of-squares (SoS) reformulations for minimizing the objective function $f$. In this method, they approximate $f$ using a multivariate polynomial model $m_c: \mathbb{R}^n \rightarrow \mathbb{R}$, referred to as the SoS Taylor model. This model comprises a $p$th-order Taylor expansion, a second-order perturbation term if the iteration is locally nonconvex (i.e., $\nabla_x^2f(x_k)$ not positive definite), and a regularization term. In locally convex iterations (i.e., $\nabla_x^2f(x_k) \succ 0$), the SoS Taylor model precisely corresponds to the AR$p$ sub-problem in \eqref{subprob}. 
Note that the regularization term for the SoS Taylor model is designed to be sufficiently large to ensure that $m_c$ is SoS-convex\footnote{The definition of SoS polynomials and SoS-convex can be found in Definitions \ref{def sos polynomial}--\ref{def sos convex}.}. The benefit of using the SoS Taylor model is that both the process of finding the $\sigma_k$ to make $m_c$ SoS convex and the process of minimizing the corresponding SoS-convex polynomial $m_c$ can be reformulated as semidefinite programmings (SDP) \cite{ahmadi2023higher, vandenberghe1996semidefinite}. 
SDP is a well-studied subclass of
convex optimization problems that can be solved to arbitrary precision within polynomial time \cite{vandenberghe1996semidefinite}.
Algorithm 1 in \cite{ahmadi2023higher} presents a minimization algorithm that employs an SoS Taylor model at each iteration. 
\cite{ahmadi2023higher} analyzed the local convergence rate and the global convergence of this method (for odd $p \ge 3$) for strongly convex objective functions. However, the global convergence and complexity for a general nonconvex objective function remain unknown.

This paper integrates the SoS Taylor model proposed in \cite{ahmadi2023higher} with the adaptive regularization techniques in \cite{birgin2017worst, cartis2020sharp, cartis2020concise, Nesterov2021implementable}. We introduce an algorithmic framework that combines the SoS Taylor model and adaptive regularization (see Algorithm \ref{arp + SoS algo}). The paper presents two key contributions.

\begin{enumerate}

\item We present both global convergence and worst-case evaluation complexity results for Algorithm \ref{arp + SoS algo}. We show that Algorithm \ref{arp + SoS algo} yields a monotonically decreasing sequence of iterates, $\{f(x_k)\}_{k \ge 0}$, and converges to an $\epsilon$-approximate first-order minimizer of $f$. The function value reduction in each strongly convex iteration is $\mathcal{O}(\epsilon^{\frac{p+1}{p}})$ for odd $p > 5$, and $\mathcal{O}(\epsilon^{\frac{p+3}{p+1}})$ for even $p > 5$. Note that we have improved function value reduction for $p = 3, 5$ at $\mathcal{O}(\epsilon^{\frac{2}{p}})$ and for $p = 4$ at $\mathcal{O}(\epsilon^{\frac{4}{5}})$. In locally nonconvex cases, the error reduction per iteration is $\mathcal{O}(\epsilon^2)$. Therefore, for a general nonconvex function, the overall evaluation complexity bound\footnote{This overall bound matches in accuracy order that for high-order trust-region methods \cite{cartis2022evaluation}, and it is not surprising to obtain this bound here, as we shall see, due to the quadratic regularisation term added in the model in nonconvex iterations.} is $\mathcal{O}(\epsilon^{-2})$.  However, for strongly convex functions $f(x)$, we establish an improved evaluation complexity bound for the $\epsilon$-approximate minimum value of convex functions as $\mathcal{O}(\epsilon^{-\frac{1}{p}})$ for odd $p$ and $\mathcal{O}(\epsilon^{-\frac{p+1}{p}})$ for even $p$.
To the best of our knowledge, for a generally nonconvex objective function, this is the first paper to offer a global convergence and a global rate result for a tensor method with a tractable sub-problem; these bounds are not of optimal or best-known accuracy order for the respective classes, but open the way for such improvements.

\item The perturbation size for locally nonconvex iterations and its effect on the complexity bound remain unspecified in \cite{ahmadi2023higher}. In this paper, we use an $\epsilon$-dependent parameter (denoted here as $\delta$) to act as the perturbation for nonconvex iterations. We provide a detailed analysis of how to choose $\delta$ as a function of $\epsilon$. Additionally, we explore how the complexity bound varies with the order of the Taylor approximation, $p$, and the magnitude of perturbation, $\delta$.
\end{enumerate}

The paper is structured as follows. In Section \ref{sec The SoS Taylor model}, we introduce the SOS Taylor model and the related adaptive regularization framework. In Section \ref{sec: upper bound for sigma}, we demonstrate the existence of an upper bound for the regularization parameter across iterations.
In Section \ref{sec: convergence}, we provide guaranteed bounds for function value reduction in successful iterations for the three convexity cases, along with the overall complexity bound for Algorithm \ref{arp + SoS algo}. We also provide results demonstrating improved function value reduction specifically for the cases of $3\le p \le 5$ and improved global convergence rate for strongly convex objective functions. In Section \ref{sec Numerics}, we present our preliminary numerical illustration of theoretical bounds, and in Section \ref{Conclusion}, we present our conclusions.

On the topic of high-order Newton methods, Gould et al. introduced the tensor-Newton approximation \cite{gould2017higher}. Convergence and evaluation-complexity analysis of a regularized tensor-Newton method for solving nonlinear least-squares problems are detailed in \cite{gould2019convergence}. Additionally, earlier work by Schnabel and Frank \cite{schnabel1991tensor, schnabel1984tensor} resulted in a practical tensor algorithm for solving nonlinear least-squares problems. The local convergence rate of AR$p$ variants is discussed in \cite{doikov2022local} which matches the local convergence result for the higher order Newton method in \cite{ahmadi2023higher}. Comprehensive analysis and literature survey on methods utilizing third derivatives and beyond are given in \cite[Ch 4]{cartis2022evaluation}.
For global optimization of polynomial functions, especially quartic polynomials, finding their global minimum is an NP-hard problem \cite{luo2010semidefinite}. However, \cite{qi2004global} presents a global descent algorithm for finding the global minimum of quartic polynomials in two variables. Another strategy involves employing branch-and-bound algorithms, as described in \cite{cartis2015branching, fowkes2013branch, horst2013global, kvasov2009univariate, neumaier2004complete}. These algorithms partition the feasible region recursively and construct nonconvex quadratic or cubic lower bounds.

\section{The SoS Taylor Model with Adaptive Regularization}
\label{sec The SoS Taylor model}

In this section, we introduce the SOS Taylor model and outline the algorithmic framework that integrates the SOS Taylor model with adaptive regularization.
In the locally strongly convex and nonconvex iterations, the SoS Taylor model is consistent with that in \cite{ahmadi2023higher}. Here, we extend the model to the locally nearly strongly convex case, which needs to be addressed when optimizing general objectives. For notational simplicity, we denote $f_k = f(x_k) \in \R$, $g_k = \nabla_x f(x_k) \in \R^n$ and $ H_k  = \nabla_x^2 f(x_k) \in \R^{ n \times n}$.

\begin{tcolorbox}[breakable, enhanced, title = The $p$th order  SoS Taylor model]
Set $p':= p+1$ if $p$ is odd and $p':= p+2$ if $p$ is even. Note that $p\ge 3$.
Let $\delta = \epsilon^a > 0$ with $ a \in [0, \frac{1}{2}]$. Note that $0 < \epsilon, \delta < 1$. 

\begin{enumerate}
\item If $H_k$ is positive definite, (i.e.  the smallest eigenvalue of Hessian satisfies $\lambda_{\min}[H_k]\ge\delta > 0$), we say that $f$ is \textbf{locally strongly convex} at $x_k$, and define the SoS Taylor model as
\begin{eqnarray}
m_c(x_k, s) =   T_{p}(x_k, s) + \frac{\sigma_k }{p'}\sigma_k \|s\|^{p'} = f_k+ g_k^Ts + \frac{1}{2}\underbrace{H_k}_{:= \bar{H}_k} [s]^2
+ \sum_{j=3}^p \frac{1}{j!} \nabla_x^j f(x_k)[s]^j +   \frac{\sigma_k }{p'}\|s\|^{p'}. 
\label{convex  SoS Taylor model}
\end{eqnarray}

\item If $H_k$ is not positive definite (i.e. $\lambda_{\min}[H_k]\le 0$), we say that $f$ is \textbf{locally nonconvex} at $x_k$,  and define the SoS Taylor model as
\begin{eqnarray} 
m_c(x_k, s) &=&    T_{p}(x_k, s)  + \frac{1}{2} \bigg(-\lambda_{\min}[H_k] +\delta\bigg) \|s\|^2+ \sigma_k \|s\|^{p+1}\notag
\\ &=&   f_k+ g_k^Ts + \frac{1}{2}\underbrace{\bigg(H_k-\lambda_{\min}[H_k] I_n+\delta I_n\bigg)  }_{:=\bar{H}_k}[s]^2
+ \sum_{j=3}^p \frac{1}{j!} \nabla_x^j f(x_k)[s]^j +   \frac{\sigma_k }{p'}\|s\|^{p'}. 
\label{nonconvex  SoS Taylor model}
\end{eqnarray}

\item If $H_k$ is nearly strongly convex, (i.e.   $0<\lambda_{\min}[H_k]\le\delta $), we say $f$ is \textbf{locally nearly strongly convex} at $x_k$,  and define the SoS Taylor model as
\begin{eqnarray}
m_c(x_k, s) &=&  T_{p}(x_k, s) + \frac{\delta}{2} \|s\|^2 + \frac{\sigma_k }{p'} \|s\|^{p'}  \notag
\\ &=&  f_k+  g_k^Ts + \frac{1}{2} \underbrace{(H_k+\delta I_n)}_{:=\bar{H}_k} [s]^2
+ \sum_{j=3}^p \frac{1}{j!} \nabla_x^j f(x_k)[s]^j  +   \frac{\sigma_k }{p'}\|s\|^{p'}. 
\label{nearly strongly convex  SoS Taylor model}
\end{eqnarray}
\end{enumerate}
\end{tcolorbox}

\begin{remark}
\label{remark bar H}
An alternative formulation of Case 3 involves setting $\bar{H}_k :=H_k-\lambda_{\min}[H_k] +\delta I_n$. Despite the perturbation term, $-\lambda_{\min}(H_k) + \delta  I_n$, becoming more complicated, this formulation offers an advantage. When $\lambda_{\min}[H_k] = \delta$, both $\bar{H}_k$ from Case 1 and $\bar{H}_k$ from Case 3 have the same smallest eigenvalue as $\delta$. 
It's noteworthy that both formulations ensure $\lambda_{\min}(\bar{H}_k) \ge \delta$. This is the crucial inequality for the proof of convergence and complexity analysis. Therefore, in the alternative formulation, the same order of convergence and complexity bounds in Section \ref{sec: convergence} would hold. For the simplicity of the perturbation term in Case 3, we will use $\bar{H}_k :=H_k +\delta I_n$  in Case 3 throughout the paper. 
\end{remark}

The key idea of the SoS Taylor model with Adaptive Regularization is to choose $\sigma_k > 0$ such that $m_c(x_k, s)$ is SoS-convex. To facilitate this, we provide definitions related to SoS, SoS-matrix, and SoS-convexity, following \cite{ahmadi2023higher,ahmadi2013complete, kojima2003sums}.

\begin{definition} \textbf{(SoS polynomial)}
    A ${\tilde{p}}$-degree polynomial $q(s): \R^n \rightarrow \R$ where $s \in \R^n$ is a sum of squares (SoS) if there exist polynomials $\Tilde{q}_1, \dotsc, \Tilde{q}_r: \R^n \rightarrow \R$, for some  $r\in \mathbb{N}$, such that $q(s) =\sum_{j=1}^r \Tilde{q}_j(s)^2$ for all $s \in \R^n$ \cite[Def. 1]{ahmadi2023higher}. 
    \label{def sos polynomial}
\end{definition}

\noindent
Definition \ref{def sos polynomial} implies that
the degree ${\tilde{p}}$ of an SoS polynomial $q$ has to be even, and that the maximum degree of each $\Tilde{q}_j$ is $\frac{\tilde{p}}{2}$. Let $\R[s]^{n\times n}$ be the real vector space of  $n \times n$ real polynomial matrices, where each entry of such a matrix is a polynomial with real coefficients.

\begin{definition} \textbf{(SoS matrix)}
A  symmetric polynomial matrix $H(s) \in \R[s]^{n\times n}$  is an SoS-matrix if there exists a polynomial matrix $V(s) \in \R[s]^{n_1\times n}$ for some $n_1 \in \mathbb{N}$, such that $H(s) = V(s)^TV(s)$ \cite[Def. 2.4]{ahmadi2013complete}. 
\label{def sos matrix}
\end{definition}

\begin{definition} \textbf{(SoS-convex)}  
A polynomial $p(s): \R^n \rightarrow \R$ is SoS-convex if its Hessian $H(s):=\nabla^2 p(s)$ is an SoS-matrix \cite[Def. 2.4]{ahmadi2013complete}.
\label{def sos convex}
\end{definition}

\noindent   
If $p(s)$ is  SoS-convex, as per Definition \ref{def sos convex}, its Hessian is an SoS-matrix which corresponds to a positive semidefinite matrix (as per Definition \ref{def sos matrix}). Consequently, if $p(s): \R^n \rightarrow \R$ is SoS-convex, it implies that $p(s)$ is a convex function \cite{ahmadi2013complete, kojima2003sums}. 


\begin{tcolorbox}[breakable, enhanced, title = Choosing $\sigma_k$ in the  SoS Taylor model]
In all three cases, according to \cite{ahmadi2023higher}, $\sigma_k>0 $ is chosen as follows according to the following sum of squares
program. 
\begin{eqnarray}
   \bar{\sigma}_k :=\min \sigma_k, \qquad \text{s.t.} \qquad m_c(x_k, s) \qquad \text{SoS-convex.}
   \label{sos program}
\end{eqnarray}
This problem can be reformulated as a semidefinite program (SDP). For fixed $p$, the size of this SDP grows polynomially in $n$. \cite[Lemma 3 and Thm 3]{ahmadi2023higher} and Theorem \ref{sigma SoS bound} here ensure that such $\sigma_k>0 $ exists and is bounded above.
\end{tcolorbox}

\begin{remark} \textbf{(Role of $\delta$)}: In \cite{ahmadi2023higher}, the authors set $\delta = \epsilon$. Here, we introduce more flexibility by allowing $\delta$ to be chosen. Specifically, we set $\delta = \epsilon^a > 0$ with $ a \in [0, \frac{1}{2}]$. We analyze this choice in later sections (Section \ref{sec: convergence}, Theorem \ref{thm: Function Value Decrease in successful iter}, and Remark \ref{choice of delta}) and discuss how to select $\delta$ as a function of $\epsilon$.  A discussion for the bound $a \ge \frac{1}{2}$ is provided in Remark \ref{remark  a>1/2} and Appendix \ref{appendix a>1/2}, respectively.
\end{remark}


\begin{remark} \textbf{(Power of regularization)}
In the AR$p$ framework \cite{birgin2017worst, carmon2020lower, cartis2022evaluation, cartis2020sharp, cartis2020concise}, the Taylor model is typically regularized by  $p'=p+1$ order. The decision to use a $p'=p+2$ order of regularization when $p$ is even is motivated by the properties of the SOS Taylor model. As indicated by \cite{ahmadi2023higher}, even power in the regularization term is necessary to ensure that the model $m_p$ can be SOS-convex.
\end{remark}

The minimization algorithm is described in Algorithm \ref{arp + SoS algo}.  
Algorithm \ref{arp + SoS algo} uses the SoS Taylor model as constructed in \eqref{convex  SoS Taylor model}, \eqref{nonconvex  SoS Taylor model} or \eqref{nearly strongly convex  SoS Taylor model} depending on the convexity case\footnote{For more on the choice of $a$ which dictates the convexity case, see  Remark \ref{choice of delta} and Remark \ref{special choice of delta}.}
together with an adaptive regularization framework.

\begin{algorithm}[H]
\caption{Adaptive $p$th-order Regularization Method with SoS Taylor Models (AR$p$+SoS)} 
\textit{Initialization}: 
Initialize $x_0 = \boldsymbol{0} \in \R^n$, $k=0$, $\sigma^{r}_k = 0$, $\sigma_{\min}>0$,  $\sigma^{r}_0 = \sigma_{\min}$,  $\eta > 0$, $\gamma_1 > 1 > \gamma_2 >0$. 

\textit{Input}: Provided an order $p \ge 3$, 
set $p':= p+1$ for odd $p$ and  $p':= p+2$ for even $p$. 
\\ Derivatives at $x_0 = \boldsymbol{0}$, $\{f(x_0), \nabla_x f(x_0),\nabla_x^2 f(x_0), \dotsc, \nabla_x^p f(x_0)\}$.   
\\ An accuracy level $\epsilon$ such that $0 < \epsilon < 1$ and $\delta = \epsilon^a > 0$. 

\textbf{Step 1: Test for termination.} If $\|\nabla_x f(x_k)\| \le \epsilon$, terminate with the $\epsilon$-approximate first-order minimizer $s_\epsilon = s_k$.

\textbf{Step 2: Step computation.}
\\Compute $\lambda_{\min}[\nabla_x^2 f(x_k)]$ to determine the convexity case.  Construct $m_c$ as \eqref{convex  SoS Taylor model}, \eqref{nonconvex  SoS Taylor model} or \eqref{nearly strongly convex  SoS Taylor model} depending on the convexity case. 
Compute $\bar{\sigma}_k$ using \eqref{sos program}. Set $\sigma_k = \max\{\bar{\sigma}_k, \sigma^{r}_k\}. $  
\\Compute $s_k = \argmin_{s \in \R^n} m_c(x_k, s)$. 

\textbf{Step 3: Acceptance of trial point.} Compute $T_p(x_k, s_k)$ and define
\begin{equation}
\tag{Ratio Test}
\rho_k = \frac{f(x_k) - f(x_k+s_k)}{f(x_k) - T_p(x_k,s_k)}.
\label{ratio test}
\end{equation}

\textbf{Step 4: Regularization parameter update.} 
\eIf {$\rho_k >\eta$ } 
{Successful iteration. $x_{k+1} := x_k+s_k$, decrease regularization $\sigma^{r}_{k+1}=\max\{ {\gamma_2}\sigma_k, \sigma_{\min}\}$.}
{Unsuccessful iteration. $x_{k+1} := x_k$, increase regularization $\sigma^{r}_{k+1}= \gamma_1 \sigma_k$.}
$k:=k+1$
\label{arp + SoS algo}
\end{algorithm}

\begin{remark} {(\textbf{Minimization of $m_c$}) 
In Algorithm \ref{arp + SoS algo}, we assume that we find a (global) minimizer of $m_c$, such that $s_k = \argmin_{s \in \R^n} m_c(x_k, s)$. 
Since $m_c$ is strongly convex and bounded below by construction, the existence and uniqueness of the minimizer are guaranteed. Therefore, any first-order minimizer is also the global minimum. The minimizer can be found by solving a semidefinite program (SDP) of size polynomial in $n$ \cite[Thm 2]{ahmadi2023higher}  and can be solved to arbitrary precision within polynomial time. 
While our proof (in particular Lemma \ref{Lemma Lower Bound for Step Size}) assumes exact first-order minimization of $\nabla m_c(x_k, s) = 0$, similar bounds can be derived when considering approximate minimization of $m_c$ with a step-termination condition, such as $\|\nabla m_c(x_k, s_k)\| \le \theta \|s_k\|^{p'-1}$ for some $\theta \in (0, 1)$. Such step-termination conditions are also used in \cite{birgin2017use, birgin2017worst,cartis2022evaluation, cartis2020concise}, with alternative variants available \cite{cartis2011adaptive}.}
\label{remark solver s_c}
\end{remark}

\section{An Upper Bound on the Regularization Parameter}
\label{sec: upper bound for sigma}

In this section, we show that there exists an iteration independent uniform upper bound for $\sigma_k = \max\{\bar{\sigma}_k, \sigma^{r}_k\} $. Note that the parameter $\bar{\sigma}_k$ that is computed in \eqref{sos program} can be considered as the smallest possible value such that $m_c(x_k, s)$ is SoS-convex, making the sub-problem tractable, while $\sigma^{r}_k$ is the regularization parameter generated in Step 4 of the adaptive regularization framework (Algorithm \ref{arp + SoS algo}) that monitors sufficient decrease and overall progress of the algorithm towards convergence.
To prove this result, we first introduce the following definitions, assumptions, and technical Lemmas.

\begin{assumption}
$f\in \C^{p, 1}(\R^n)$ with $p \ge 3$ which means that $f$ is $p$-times continuously differentiable, bounded below by a constant $f_\text{low}$ and the $p$th derivative of $f$ is globally Lipschitz continuous. Namely, there exists a constant $L > 0$ such that, for all
$x, y \in \R^n$, $\|\nabla_x^p f(x) - \nabla_x^p f(y)\| \le (p-1)! L \|x-y\|.$
\label{assumption Liptz}
\end{assumption}

\noindent Note that Assumption \ref{assumption Liptz} is a standard assumption in the analysis of algorithms that utilize high-order derivative information \cite{birgin2017worst, cartis2007adaptive, cartis2011adaptive, cartis2020concise}.  

\begin{definition} Let $ T \in \R^{n^j}$ be a $j$th-order tensor. Following \cite{cartis2020concise}, the tensor norm of $ T$ is
$$\|T\|_{[j]} := \max_{\|v_1\|=\dotsc=\|v_j\|=1} \big|T[v_1, \dotsc, v_j]\big|$$ 
where $\| \cdot \|$ is the Euclidean norm and $v_1, \dotsc, v_j \in \R^n$ are vectors. Note that $ T[v_1, \dotsc, v_j]$ stands for applying the $j$th-order tensor $T$ to the vectors $v_1, \dotsc, v_j$. 
\label{tensor norm}
\end{definition}

\begin{assumption}
We assume that the tensor norms $ \big\|\nabla_x^j f(x_k) \big\|_{[j]}$ are uniformly bounded at all iterates $x_k$, such that
$$
 \big\|\nabla_x^j f(x_k) \big\|_{[j]}:=  \max_{\|v_1\|=\dotsc=\|v_j\|=1} \bigg| \nabla_x^j f(x_k) [v_1,\dotsc,v_j]  \bigg|\le \Lambda_j,
$$
for all $k \ge 0$ and $j = 1, \dotsc, p$. 
\label{assumption bounded hessian}
\end{assumption}

\noindent Note that Assumption \ref{assumption bounded hessian} is only needed for iterates $x_k$ produced by the algorithm. 

\begin{definition}
Let $q(s): \R^n \rightarrow \R$ be a polynomial. Following \cite{ahmadi2022complexity}, we define $\|q(s)\|_{\infty^*}$ as the largest absolute value of the coefficients of $q(s)$ in the standard monomial basis. 
\label{defition q*}
\end{definition}

\begin{lemma} \textbf{(Related norms)}
Under Assumption \ref{assumption bounded hessian}, the following statements hold. 
\begin{enumerate}
    \item  The coefficients of $\frac{1}{j!}\nabla_x^j f(x_k) [s]^j$ in the standard monomial basis  at the iterates $x_k$ are bounded, such that $\big\|\frac{1}{j!}\nabla_x^j f(x_k) [s]^j\big\|_{\infty^*}  \le \Lambda_j$. 
    \item The leftmost eigenvalue of the second-order derivatives of $f$ at the iterates $x_k$ is bounded, such that  $\max_{k}\big|\lambda_{\min}[\nabla_x^2 f(x_k)]\big|   \le \Lambda_2.$
\end{enumerate}
\label{lemma norm}
\end{lemma}

\begin{proof}
    The proof of Lemma \ref{lemma norm} can be found in Appendix \ref{proof of lemma norm}
\end{proof}

\noindent
The subsequent Lemma and assumption can be used to prove that $\sigma_k$ generated by the SoS-convex paradigm remains bounded regardless of the number of iterations.

\begin{lemma}(\cite[Lemma 2]{ahmadi2023higher}) Let $q: \R^n \rightarrow \R$ be a multivariate polynomial expressed in the standard monomial basis and $\|\cdot\|_{\infty^*}$ be defined as in Definition \ref{defition q*}. There exists a constant $R>0$, such that for any polynomial $q(s)$ of degree at most $p$ with $\|q\|_{\infty^*} \le R$, the polynomial $\|s\|^2 + \|s\|^{p'} + q(s)$ is SoS-convex, where $p':= p+1$ if $p$ is odd and $p':= p+2$ if $p$ is even. 
\label{lemma ali}
\end{lemma}

\begin{remark}
\cite[Lemma 2]{ahmadi2023higher} gives that for an even integer $p'$, $\|s\|^2 + \|s\|^{p'}$ lies within the interior of the cone of sos-convex polynomials in $n$ variables and of degree at most $p'$. This result implies Lemma \ref{lemma ali}. Importantly, $R$ in Lemma \ref{lemma ali} is a constant that is independent of the iteration $k$, $\epsilon$, and $\delta$.
\end{remark}

To demonstrate that $\sigma_k = \max\{\bar{\sigma}_k, \sigma^{r}_k\}$ is uniformly bounded above regardless of the iteration $k$, we outline our proof in two steps, with some numerical illustration on the size of $\sigma_k$ are presented in Section \ref{sec Numerics}.

\begin{enumerate}
    \item In Section \ref{sec: uniform upper bound 1}, we first show that there exists a uniform upper bound for $\bar{\sigma}_k$ such that $m_c(x_k, s)$ is SoS-convex. Note that in \cite{ahmadi2023higher}, only the strongly convex case and the locally nonconvex case are discussed. We extend this result to include the locally nearly strongly convex case, so that all possible cases are now addressed. We also derive an explicit expression for the upper bound on $\bar{\sigma}_k$ (Corollary \ref{bound for sigma}). 
    \item Then, in Section \ref{sec: uniform upper bound 2}, we show that there exists a uniform upper bound on $\sigma^{r}_k$ to ensure that the next iteration gives a sufficient decrease in $f$ (Theorem \ref{thm general bound for sigma odd} for odd $p$ and Theorem \ref{thm general bound for sigma odd} for even $p$).
\end{enumerate}

\subsection{A Uniform Upper Bound on $\bar{\sigma}_k$ }
\label{sec: uniform upper bound 1}

In this subsection, we establish a uniform upper bound on $\bar{\sigma}_k$. 

\begin{theorem} (\cite[Theorem 3]{ahmadi2023higher}\footnote{In \cite{ahmadi2023higher}, slightly different notations are employed, where the cubic and higher terms are defined as $C(s):= \frac{2}{\delta}\sum_{j=3}^p \frac{1}{j!} \nabla_x^j f(x_k)[s]^j$ and $M = (2/\delta) M_k$.})
Let $C_{k,p}(s):= \sum_{j=3}^p \frac{1}{j!} \nabla_x^j f(x_k)[s]^j$ and let  $M_k :=\|C_{k,p}(s)\|_{\infty^*}$ where  $\|\cdot\|_{\infty^*}$ is defined in Definition \ref{defition q*}. Let $f\in \C^{p, 1}(\R^n)$ and $\bar{\sigma}_k$ is computed using \eqref{sos program} in Algorithm \ref{arp + SoS algo}. Then, 
\begin{eqnarray}
   \frac{\bar{\sigma}_k}{p'}\le \frac{\delta}{2} \min\bigg\{1, \frac{2 M_k}{\delta R}\bigg\}^{p'-2} 
   \label{ali bound for sigma}
\end{eqnarray}
where $R$ is defined in Lemma \ref{lemma ali} and is a constant independent of the iteration $k$, $\epsilon$, and $\delta$.
\label{thm ali}
\end{theorem}

\begin{corollary}
\label{bound for sigma} 
\textbf{(Upper bound for $ \bar{\sigma}_k$)}
Let $f\in \C^{p, 1}(\R^n)$ and $\bar{\sigma}_k$ is computed using \eqref{sos program} in Algorithm \ref{arp + SoS algo}. Under Assumption \ref{assumption bounded hessian}, in all convexity cases and for all $k \ge 0$, 
\begin{eqnarray}
    \bar{\sigma}_k \le  \max \bigg\{\frac{\delta  p'}{2}, \quad C_R     \delta^{3-p'} \bigg\}, \quad {\rm where}\quad  
    C_R = \frac{p'}{2}   \bigg(\sum_{j = 3}^p \frac{2 \Lambda_j}{R}\bigg)^{p'-2}. 
    \label{sigma SoS bound}
\end{eqnarray}
Note that $C_R$ is independent of the iteration count $k$, the tolerance $\epsilon$, and the perturbation $\delta$. 
\end{corollary}

\begin{proof}
Using  Assumption \ref{assumption bounded hessian} and Lemma \ref{lemma norm}, we deduce that 
\begin{eqnarray} 
M_k = \|C_{k,p}(s)\|_{\infty^*} \le \sum_{j = 3}^p \bigg\|\frac{1}{j!}\nabla_x^j f(x_k) [s]^j\bigg\|_{\infty^*} \le  \sum_{j = 3}^p \Lambda_j.
\label{our M}
\end{eqnarray}
In the locally strongly convex case, we have $\bar{H}_k: =H_k \succeq \delta I_n $. Similarly, in the locally nonconvex case, we have $\bar{H}_k:= H_k-\lambda_{\min}[H_k] I_n+\delta I_n  \succeq \delta I_n $, and in the locally nearly strongly convex case, we have $\bar{H}_k:= H_k +\delta I_n \succeq \delta I_n$. Therefore, we can apply Theorem \ref{thm ali} in all three cases. Substituting \eqref{our M} into \eqref{ali bound for sigma}, we obtain that
$$
    \frac{\bar{\sigma}_k}{p'}\le \frac{\delta}{2} \min\bigg\{1, \frac{(2/\delta) M_k}{R}\bigg\}^{p'-2} 
    \le \max \bigg\{ \frac{\delta}{2}, \quad \bigg(\sum_{j = 3}^p \frac{\Lambda_j}{R} \bigg)^{p'-2}  \bigg(\frac{\delta}{2}\bigg)^{3-p'}\bigg\}.
$$
Setting  $C_R$ as stated in the lemma, we deduce the required bound.
\end{proof}

\subsection{A Uniform Upper Bound on $\sigma^r_k$ }
\label{sec: uniform upper bound 2}

 We first introduce two technical lemmas.

\begin{lemma}\cite[Lemma 2.1]{cartis2020concise} 
   Let $f\in \C^{p, 1}(\R^n)$ and $T_p(x, s)$ be the $p$th Taylor approximation of $f (x + s)$ about $x$. Then, under Assumption \ref{assumption Liptz}, for all $x, s \in \R^n$, 
\begin{eqnarray}
      |f(x+s) - T_p(x, s) | &\le&  \frac{L}{p+1}\|s\|^{p+1},
            \label{lip bound 1}
      \\ \big\|\nabla_x f(x+s) - \nabla_s T_p(x, s) \big\|_{[p]}  &\le& L \|s\|^{p},
      \label{lip bound 2}
\end{eqnarray}
where $ \| \cdot \|_{[p]} $ is defined as in Definition \ref{tensor norm}.
\label{tech lemma 1}
\end{lemma}
\begin{proof}
  The proof can be found in Lemma 2.1 and Appendix A.1  \cite{cartis2020concise}.
\end{proof}

\begin{lemma}  \label{tech lemma 2}  Let $f\in \C^{p, 1}(\R^n)$.
 The mechanism of Algorithm \ref{arp + SoS algo} guarantees that, for all $k \ge 0$,
\begin{eqnarray}
&\text{1: locally strongly convex:} &  f(x_k) - T_p(x_k, s_k) \ge \frac{\sigma_k}{p'} \|s_k\|^{p'},
      \label{temp mT 1}
\\&\text{2: Locally nonconvex:}&  f(x_k) - T_p(x_k, s_k) \ge \big(-\lambda_{\min}[H_k]+\delta\big)\|s_k\|^{2}+  \frac{\sigma_k}{p'} \|s_k\|^{p'},
      \label{temp mT 2}
\\&\text{3: Locally nearly strongly convex:}&  f(x_k) - T_p(x_k, s_k) \ge \delta \|s_k\|^{2}+ \frac{\sigma_k}{p'} \|s_k\|^{p'} .
      \label{temp mT 3}
\end{eqnarray}
\end{lemma}
\begin{proof} Since  $s_k = \argmin_{s \in \R^n} m_c(x_k, s)$, we have  $m_c(0) - m_c(x_k, s_k) \ge 0$. Thus, 
$$
   f(x_k) - T_p(x_k, s_k) = \underbrace{\big[f(x_k) - m_c(x_k, s_k) \big]}_{\ge 0} + \big[m_c(x_k, s_k)- T_p(x_k, s_k)\big] \ge m_c(x_k, s_k)- T_p(x_k, s_k).
$$
\begin{itemize}
\item In locally strongly convex cases, according \eqref{convex  SoS Taylor model},  
$
    m_c(x_k, s_k)- T_p(x_k, s_k) = \frac{\sigma_k}{p'} \|s_k\|^{p'}.
$
\item  In locally nonconvex cases, according \eqref{nonconvex  SoS Taylor model},  
$
m_c(x_k, s_k)- T_p(x_k, s_k) =  \big(-\lambda_{\min}[H_k]+\delta\big)\|s_k\|^{2}+ \frac{\sigma_k}{p'} \|s_k\|^{p'}.
$
\item  In locally nearly strongly convex cases, according to \eqref{nearly strongly convex  SoS Taylor model},  
$
m_c(x_k, s_k)- T_p(x_k, s_k) =  \delta\|s_k\|^{2}+ \frac{\sigma_k}{p'} \|s_k\|^{p'}.
$
\end{itemize}
\end{proof}

\begin{lemma} 
\label{lemma: success sigma k}
\textbf{(Upper bound for $\sigma_k$)} 
Let $f\in \C^{p, 1}(\R^n)$. Under Assumption \ref{assumption Liptz}, if
\begin{eqnarray}
    \sigma_k \ge \max \big\{\sigma_0, \quad  A_0\|s_k\|^{-\iota} \big\},
    \label{sigma bound L}
\end{eqnarray}
where $A_0 = \frac{L p'}{(p+1)(1-\eta)}$, $\iota=0$ for odd $p$, and $\iota=1$ for even $p$, then, $\rho_k \ge \eta$ and \eqref{ratio test} for iteration $k$ is successful.
Note that\footnote{\eqref{sigma bound L} matches the bound in the universal ARP$r$ framework (\cite[Lemma 3.4]{cartis2019universal}) by setting the parameter in  \cite{cartis2019universal} as $\beta_p = 1$ and $r = p+1$ for odd $p$ and $r = p+2$ for odd $p$. } if $p$ is odd,  \eqref{sigma bound L} becomes $ \sigma_k \ge \max \big\{\sigma_0, \frac{\gamma_1 L}{(1-\eta)} \big\} $; if $p$ is even, \eqref{sigma bound L} becomes  $ \sigma_k \ge \max \big\{\sigma_0, \frac{\gamma_1 L(p+2)}{(p+1)(1-\eta)\|s_k\|} \big\}. $ 
\end{lemma}

\begin{proof} The proof follows a similar setup as \cite[Lemma 3.2]{cartis2020concise},  \cite[Lemma 3.4]{cartis2019universal} and  \cite[Lemma 2.2]{birgin2017worst}. From \eqref{temp mT 1}--\eqref{temp mT 3}, we deduce that in all three cases $ f(x_k) - T_p(x_k, s_k) \ge \frac{\sigma_k}{p'} \|s_k\|^{p'}.$ Also, using \eqref{lip bound 1}, we may then deduce that
$$
|\rho_k-1| \le \frac{|f(x_k+s_k) - T_p(x_k, s_k)|}{|f(x_k) - T_p(x_k, s_k)|} \le \frac{ Lp'\|s_k\|^{p+1}}{(p+1) \sigma_k \|s_k\|^{p'}} \le \frac{L p' }{(p+1) \sigma_k } \|s_k\|^{p+1-p'}. 
$$
Let  $\iota :=p'- (p+1)$ and  $\iota=0$ for odd $p$, and $\iota=1$ for even $p$.  Assume $\sigma_k \ge \frac{L p'}{(p+1)(1-\eta)} \|s_k\|^{-\iota}$,  we deduce that $|\rho_k-1|  \le 1 - \eta$ and thus $\rho_k \ge \eta$. Then, iteration $k$ is successful.
\end{proof}

\begin{tcolorbox}
\begin{theorem} \textbf{(Uniform upper bound for $\sigma_k$ when $p$ is odd.)}
Let $f\in \C^{p, 1}(\R^n)$ and Assumptions \ref{assumption Liptz}--\ref{assumption bounded hessian} hold. If $p \ge 3$ is odd and $p'= p+1$,  for all $k \ge 0$, 
\begin{eqnarray}
    \sigma_k \le \max \big\{\sigma_0, \quad \gamma_1 A_0, \quad p'\frac{\delta}{2},  \quad  C_R   \delta^{3-p'}   \big\}. 
    \label{general bound for sigma odd}
\end{eqnarray}
where\footnote{In the general AR$p$ framework \cite{cartis2020sharp, cartis2020concise}, we have $p'= p+1$. Also, we do not need the upper bound in \eqref{sigma SoS bound}, and so
$    \sigma_k \le  \max \big\{\sigma_0, \frac{\gamma_1 L}{(1-\eta)}\big\}. $} $A_0:= \frac{L p'}{(p+1)(1-\eta)}$. 
$\sigma_0, \gamma_1, L, p, \eta, C_R, A_0$ are independent of the tolerance  $\epsilon$, and  $\delta = \epsilon^a$ with $a \in [0, \frac{1}{2}]$. Moreover, for any $0  < \epsilon < 1$, we have
\begin{eqnarray}
   \sigma_k \le   \hat{C}_R \epsilon^{a(2-p)} =: \sigma_{\max}
    \label{general bound for sigma odd 2}
\end{eqnarray}
where $\hat{C}_R =  \max \big\{\sigma_0, \gamma_1 A_0, p'\frac{\delta}{2},   C_R  \big\}$. Note that $ C_R = \frac{p'}{2} (\sum_{j = 3}^p \frac{2 \Lambda_j}{R})^{p'-2}$ as defined in \eqref{sigma SoS bound}. 
\label{thm general bound for sigma odd}
\end{theorem}
\end{tcolorbox}

\begin{proof}
Assume that $p \ge 3$ is odd and $p'= p+1$. To prove the first result, assume $\sigma_k \ge \frac{\gamma_1 L p'}{(p+1)(1-\eta)}$,  according to Lemma \ref{lemma: success sigma k}, we deduce that iteration $k$ is successful and $\sigma^r_{k+1} \le \sigma_k$. As a consequence, the mechanism of the algorithm ensures that 
$$  
\sigma^r_k \le \max \bigg\{\sigma_0, \frac{L p' \gamma_1}{(p+1)(1-\eta)} \bigg\}
$$
for all $k$. Combining with the bound \eqref{sigma SoS bound} in Corollary \ref{bound for sigma}, we have $\sigma_k = \max\{\bar{\sigma}_k, \sigma^{r}_k\}$ which gives \eqref{general bound for sigma odd}.
To prove the second result, using \eqref{general bound for sigma odd} and $\hat{C}_R =  \max \big\{\sigma_0, \gamma_1 A_0, p'\frac{\delta}{2},   C_R  \big\}$, we have  
\begin{eqnarray}
\sigma_k \le  \hat{C}_R \max \big\{1,  \delta^{3-p'} \} =  \hat{C}_R \max \big\{1,\epsilon^{(2-p)a} \}.
\label{odd bound a}
\end{eqnarray}
By considering  $0 < \epsilon < 1$, we obtain  \eqref{general bound for sigma odd 2}.
\end{proof}

\begin{lemma}Let $f\in \C^{p, 1}(\R^n)$ and Assumptions \ref{assumption Liptz}--\ref{assumption bounded hessian} hold.  If $p \ge 4$ is even and $p'= p+2$, for all $k \ge 0$ and any $0< \epsilon <1$,  
    \begin{eqnarray}
    \sigma_k \ge \max\{\sigma_0, B\epsilon^{-1}\},   
    \label{bound for odd}
\end{eqnarray}
where $B = 4A_0 \max \big\{1, (\Lambda_2 + 1), (L+A_0)^{\frac{1}{p}} \}$,  then \eqref{sigma bound L} holds,  $\rho_k \ge \eta$ and \eqref{ratio test} for iteration $k$ is successful.
\label{lemma: success sigma k 2}
\end{lemma}

\begin{proof}
    We will prove our result by contradiction. We assume \eqref{sigma bound L} does not hold on iteration $k$, such that
    \begin{eqnarray}
    \sigma_k \|s_k\|< \frac{L p'}{(p+1)(1-\eta)} = A_0. 
    \label{contradiction argument}
\end{eqnarray}
    Note that while Algorithm \ref{arp + SoS algo} does not terminate, we have
    \begin{eqnarray}
\epsilon &<& \big\|\nabla_x f(x+s_k)\big\| = \big\|\nabla_x f(x+s_k) - \nabla_s m_c(x_k, s_k)  \big\| \notag
\\ &\le & \big\|\nabla_x f(x+s_k) - \nabla_s T_p(x, s_k) \big\| + \big\|\nabla_s m_c(x_k, s_k)  \ - \nabla_s T_p(x, s_k) \big\| \notag
\\ &\underset{\eqref{lip bound 2}}{\le} &  L \|s_k\|^p + \big\|\nabla_s m_c(x_k, s_k)  \ - \nabla_s T_p(x, s_k) \big\|. 
\label{step size ineq 1}
\end{eqnarray}
Using  \eqref{temp mT 1}--\eqref{temp mT 3} and \eqref{contradiction argument}, \eqref{step size ineq 1} becomes\footnote{In the universal ARP$r$ framework \cite{cartis2019universal}, the framework uses the Taylor model with no perturbation term, i.e., $\delta = 0$. Therefore, \eqref{case 1 comp bound}--\eqref{case 3 comp bound} become $\epsilon < (L+A_0) \|s_k\|^p$ and we have $\sigma_k \ge \mathcal{O}(\epsilon^{-\frac{1}{p}})$ \cite[Lemma 3.4]{cartis2019universal}.}
\begin{eqnarray}
\text{Case 1: }  &&\epsilon <L \|s_k\|^p + \sigma_k \|s_k\| \|s_k\|^p  \underset{\eqref{contradiction argument}}{\le} (L+A_0) \|s_k\|^p,
\label{case 1 comp bound}
\\\text{Case 2: }   && 
\epsilon <  \big(-\lambda_{\min}[H_k] + \delta\big)\|s_k\|+ L \|s_k\|^p + \sigma_k \|s_k\| \|s_k\|^p  \le  \big(\Lambda_2+1\big)\|s_k\|+ (L+A_0) \|s_k\|^p,
\label{case 2 comp bound}
\\\text{Case 3: }  &&
\epsilon <   \delta \|s_k\|+ L \|s_k\|^p + \sigma_k \|s_k\| \|s_k\|^p  \le  \epsilon^a\|s_k\| + (L+A_0) \|s_k\|^p,
\label{case 3 comp bound}
\end{eqnarray} 
where the last inequality in Case 2 is obtained by Lemma \ref{lemma norm} and $0< \delta = \epsilon^a < 1$. From \eqref{case 1 comp bound}, we deduce that $\|s_k\| > \big(\frac{\epsilon}{L+A_0}\big)^{\frac{1}{p}} $. From \eqref{case 2 comp bound}, we deduce that either $\frac{\epsilon}{2} <  (L+A_0) \|s_k\|^p$ or $\frac{\epsilon}{2} <   \big(\Lambda_2+1\big)\|s_k\|$. Therefore, $\|s_k\| > \frac{1}{2} \min \big\{ \frac{\epsilon}{\Lambda_2 + 1} ,  \big(\frac{\epsilon}{L+A_0}\big)^{\frac{1}{p}} \big\}. $ Similarly, we deduce from \eqref{case 3 comp bound} that $    \|s_k\| > \frac{1}{2} \min \big\{ \epsilon^{1-a},  \big(\frac{\epsilon}{L+A_0}\big)^{\frac{1}{p}} \big\}. $ Combining these bounds, we obtain that in all iterations, $\|s_k\|$ has the following lower bound,
  \begin{eqnarray*}
      \|s_k\| > \frac{1}{2} \min \bigg\{ \epsilon^{1-a}, \quad \frac{\epsilon}{\Lambda_2 + 1} , \quad \bigg(\frac{\epsilon}{L+A_0}\bigg)^{\frac{1}{p}} \bigg\} \ge B_1  \min  \big\{ \epsilon^{1-a}, \epsilon, \epsilon^{\frac{1}{p}}\} = B_1 \epsilon . 
      \label{temp 5}
  \end{eqnarray*}
  where $B_1 = \frac{1}{2}  \min \big\{1, (\Lambda_2 + 1)^{-1}, (L+A_0)^{-\frac{1}{p}} \}$ is a constant. The last inequality comes from $\epsilon \le  \epsilon^{1-a}$ and  $\epsilon \le  \epsilon^{\frac{1}{p}}$  for $0< \epsilon < 1$ and $a >0$. On the other hand, by \eqref{contradiction argument}, we have also $ \|s_k\| \le \frac{A_0}{ \sigma_k} . $ Therefore,
  \begin{eqnarray*}
     \frac{A_0}{ \sigma_k} >   B_1 \epsilon  \qquad \Rightarrow \qquad \sigma_k < 2A_0 B_1^{-1} \epsilon^{-1}
  \end{eqnarray*}
which contradicts \eqref{bound for odd}. Thus, \eqref{sigma bound L} must hold and Lemma \ref{lemma: success sigma k} implies that $k$ is successful.
\end{proof}

\begin{tcolorbox}
\begin{theorem} \textbf{(Upper bound for $\sigma_k$ when $p$ is even.)}
Let $f\in \C^{p, 1}(\R^n)$ and Assumptions \ref{assumption Liptz}--\ref{assumption bounded hessian} hold. If $p \ge 4$ is even and $p'= p+2$, for all $k \ge 0$, 
\begin{eqnarray}
    \sigma_k \le \max \big\{\sigma_0, \quad \gamma_1 B\epsilon^{-1}, \quad p'\frac{\delta}{2},  \quad  C_R   \delta^{3-p'}   \big\} 
    \label{general bound for sigma even}
\end{eqnarray}
where $\gamma_1, L, p, \eta, C_R, B$ are independent of the tolerance  $\epsilon$, and  $\delta = \epsilon^a$ with $a \in [0, \frac{1}{2}]$. Moreover, for any $0< \epsilon < 1$, we have 
\begin{eqnarray}
   \sigma_k \le \tilde{C}_R \max\big\{ \epsilon^{-1} ,   \epsilon^{a(1-p)}\big\} = \tilde{C}_R  \epsilon^{-\max\{1,   a(1-p)\} } =: \sigma_{\max}
   \label{general bound for sigma even 2}
\end{eqnarray}
where $\tilde{C}_R = \max\{  \sigma_0,  \gamma_1 B,   p'\frac{\delta}{2}, C_R \}$. Note that if $a \in \big[0, \frac{1}{p-1}]$, then  $   \sigma_k \le \tilde{C}_R \epsilon^{-1}$. If $a \in \big[\frac{1}{p-1}, \frac{1}{2}\big]$, then  $   \sigma_k \le \tilde{C}_R \epsilon^{a(1-p)}$. 
\label{thm general bound for sigma even}
\end{theorem}
\end{tcolorbox}

\begin{proof}
Assume that $p \ge 4$ is even and $p'= p+2$. To prove the first result, assume $\sigma_k \ge  \gamma_1 B\epsilon^{-1}$,  according to Lemma \ref{lemma: success sigma k 2}, we deduce that iteration $k$ is successful and $\sigma^r_{k+1} \le \sigma_k$. As a consequence, the mechanism of the algorithm ensures that 
$  
\sigma^r_k \le  \gamma_1 B\epsilon^{-1}
$
for all $k$. Combining with the bound \eqref{sigma SoS bound} in Corollary \ref{bound for sigma}, we have $\sigma_k = \max\{\bar{\sigma}_k, \sigma^{r}_k\}$ which gives \eqref{general bound for sigma even}.
To prove the second result, we start from \eqref{general bound for sigma even}, which gives us $  \sigma_k \le  \max \big\{\sigma_0,  \gamma_1 B\epsilon^{-1},  p'\frac{\delta}{2},  C_R   \epsilon^{(3-p')a} \}$.  Let  $\tilde{C}_R = \max\{  \sigma_0,  \gamma_1 B,   p'\frac{\delta}{2}, C_R \}$. Then, we have 
\begin{eqnarray}
\sigma_k \le  \tilde{C}_R \max \big\{1, \epsilon^{-1},  \delta^{3-p'} \} =  \tilde{C}_R \max \big\{1, \epsilon^{-1},  \epsilon^{(1-p)a} \}.
\label{even bound a}
\end{eqnarray}
By considering $0 < \epsilon < 1$, we obtain \eqref{general bound for sigma even 2}.
\end{proof}

\section{Convergence and Complexity Analysis}
\label{sec: convergence}

In this section, we prove that Algorithm \ref{arp + SoS algo} generates a monotonically decreasing sequence of function values, $\{f(x_k)\}_{k\ge 0}$, and that $\lim_{k\rightarrow \infty} \inf \|\nabla f(x_k)\| = 0$ (i.e., at least one limit points of $x_k$ is a stationary point of $f(x)$). By setting $a=0$, we show that, in strongly convex iterations, Algorithm \ref{arp + SoS algo} achieves a function value decrease that is at least of order $\mathcal{O}(\epsilon^{\frac{p+1}{p}})$ for odd $p > 5$, and $\mathcal{O}(\epsilon^{\frac{p+3}{p+1}})$ for even $p > 5$. Note that we have improved function value decrease for $p = 3, 5$ at $\mathcal{O}(\epsilon^{\frac{2}{p}})$ and for $p = 4$ at $\mathcal{O}(\epsilon^{\frac{4}{5}})$. 
In locally nonconvex cases, the objective reduction per iteration is at least $\mathcal{O}(\epsilon^2)$. When $f$ is globally strongly convex, we establish an enhanced evaluation complexity bound for the $\epsilon$-approximate minimum value of convex functions as $\mathcal{O}(\epsilon^{-\frac{1}{p}})$ for odd $p$, matching the result in \cite{ahmadi2023higher, Nesterov2021implementable}.

Our convergence and complexity analysis use standard techniques similar to those in \cite{cartis2007adaptive, cartis2020sharp, cartis2020concise}. Firstly, we establish lower bounds on model decrease and step length at each iteration. Note that ensuring a lower bound on the step is crucial in guaranteeing function value reduction, incorporating techniques from \cite[Lemma 2.3]{birgin2017worst} and \cite[Lemma 3.3]{cartis2020concise}. Then, we construct the lower bound for function value reduction using the relationship between model decrease characterized by \eqref{ratio test}.
While standard techniques offer valuable insights, our innovation lies in their application to the SoS Taylor model across various convexity scenarios. To the best of our knowledge, we have not encountered this type of global complexity results for a tractable cubic or higher-order sub-problem in nonconvex smooth optimization.

The section is structured as follows. In Section \ref{sec fun val decrease}, we present the function value reductions in successful iterations of Algorithm \ref{arp + SoS algo} for three different convexity cases. We analyze the impact of choosing $\delta$ as a function of $\epsilon$ on the bound for function decrease. Additionally, improved bounds for function value reduction for $3 \le p \le 5$ are provided in Sections \ref{sec Improved Function Value Decrease for p<5} and \ref{sec proof for p=4,5}, respectively. Subsequently, in Section \ref{sec overall Complexity Bound}, we provide the overall complexity bound for Algorithm \ref{arp + SoS algo} across three convexity scenarios counting both successful and unsuccessful iterations. Finally, Section \ref{sec overall convex Complexity Bound} outlines the improved overall complexity bound for a strongly convex objective function $f$.


\subsection{Bounding Objective Decrease in Successful Iterations}
\label{sec fun val decrease}

\begin{theorem} 
    \label{Lemma Lower Bound for Step Size}
\textbf{(Lower bound on the step size)} 
Let $f\in \C^{p, 1}(\R^n)$ with $p \ge 3$. Let $\delta = \epsilon^a$ for any $0 < \epsilon <1$ and $ a \in [0, \frac{1}{2}]$. Under Assumptions \ref{assumption Liptz}--\ref{assumption bounded hessian}, for all $k \ge 0$ and any $0< \epsilon <1$, 
\begin{eqnarray}
\text{Case 1: }  \quad  &&  \|s_k\| > \frac{1}{2}\min\bigg\{\bigg(\frac{\epsilon}{L}\bigg)^{\frac{1}{p}}, \bigg(\frac{\epsilon}{\sigma_k}\bigg)^{\frac{1}{p'-1}} \bigg\},
\label{step low bound convex}
\\\text{Case 2: }  \quad &&  
\|s_k\| >\hat{\kappa}_3  \epsilon, 
\label{step low bound nonconvex}
\\\text{Case 3: } \quad && \|s_k\| >  \hat{\kappa}_4  \epsilon^{1-a},
\label{step low bound nearly strongly convex}
\end{eqnarray}
\end{theorem}
where $\hat{\kappa}_3$ and $\hat{\kappa}_4$ are iteration independent and  $\epsilon$ independent  constants. 

\begin{proof} Using the same deduction as in \eqref{step size ineq 1}, we arrive at $
\epsilon <  L \|s_k\|^p + \big\|\nabla_s m_c(x_k, s_k)  - \nabla_s T_p(x, s_k) \big\|. 
$
Using  \eqref{temp mT 1}--\eqref{temp mT 3},\eqref{general bound for sigma odd 2}, \eqref{general bound for sigma even 2} and \eqref{contradiction argument}, we have
\begin{eqnarray}
\text{Case 1: }    &&\epsilon <L \|s_k\|^p + \sigma_k  \|s_k\|^{p'-1},
\label{convex epi}
\\\text{Case 2: }   && 
\epsilon <   \big(-\lambda_{\min}[H_k] + \delta\big)\|s_k\|+ L \|s_k\|^p + \sigma_k \|s_k\|^{p'-1},
\label{nonconvex epi}
\\\text{Case 3: } &&
\epsilon <    \epsilon^a \|s_k\|+ L \|s_k\|^p + \sigma_k \|s_k\|^{p'-1}.
\label{nearlyconvex epi}
\end{eqnarray} 

\noindent
\textbf{Case 1 (locally strongly convex):} Using \eqref{convex epi}, we deduce that
$
\|s_k\| >  \frac{1}{2}\min\bigg\{\big(\frac{\epsilon}{L}\big)^{\frac{1}{p}}, \big(\frac{\epsilon}{\sigma_k}\big)^{\frac{1}{p'-1}} \bigg\}.
$

\noindent
\textbf{Case 2 (locally nonconvex):}  Using \eqref{nonconvex epi} and $\sigma_k \le \sigma_{\max}$,  we deduce \\$
      \|s_k\| > \frac{1}{3} \min \bigg\{ \frac{\epsilon}{-\lambda_{\min}[H_k]  + \delta} , \big(\frac{\epsilon}{L}\big)^{\frac{1}{p}} ,  \big(\frac{\epsilon}{\sigma_{\max}}\big)^{\frac{1}{p'-1}} \bigg\}. 
$ Let $\hat{\kappa}_3 :=\frac{1}{3} \min \bigg\{(\Lambda_2  + 1)^{-1},  L^{\frac{1}{p}},   \hat{C}_R^{-\frac{1}{p}}, \tilde{C}_R^{-\frac{1}{p+1}}   \bigg\} $. 
\begin{itemize}
    \item 
If $p$ is odd, $p' = p+1$, by \eqref{general bound for sigma odd 2}, $\sigma_{\max} := \hat{C}_R  \epsilon^{a(2-p)}$, then 
$$
      \|s_k\| > 
     \underbrace{\frac{1}{3} \min \bigg\{(-\lambda_{\min}[H_k]  + \delta)^{-1},  L^{\frac{1}{p}},  \hat{C}_R^{-\frac{1}{p}}  \bigg\}}_{\ge \hat{\kappa}_3} \min \bigg\{\epsilon, \epsilon^{\frac{1}{p}} ,  \epsilon^{\frac{{a(p-2)+1}}{p}} \bigg\} \ge \hat{\kappa}_3 \epsilon 
$$
where the last inequality comes from $1 \ge  \frac{{a(p-2)+1}}{p} \ge \frac{1}{p}$ for all $p \ge 3$ and $a \in [0, \frac{1}{2}]$ and  $\min \big\{\epsilon, \epsilon^{\frac{1}{p}} ,  \epsilon^{\frac{{a(p-2)+1}}{p}} \big\}=\epsilon$  for $0<\epsilon<1$. 
\item If $p$ is even, $p' = p+2$, by \eqref{general bound for sigma even 2}, $\sigma_{\max} :=\tilde{C}_R \max\big\{ \epsilon^{-1} ,   \epsilon^{a(1-p)}\big\} $, then 
$$
      \|s_k\| >     \underbrace{\frac{1}{3} \min \bigg\{(-\lambda_{\min}[H_k]  + \delta)^{-1},  L^{\frac{1}{p}}, \tilde{C}_R^{-\frac{1}{p+1}}   \bigg\}}_{\ge \hat{\kappa}_3}  \min \bigg\{ \epsilon ,  \epsilon^{\frac{1}{p}} ,  \epsilon^{\frac{{\max\{a(p-1), 1\}+1}}{p+1}} \bigg\} \ge \hat{\kappa}_3 \epsilon
$$
where the last inequality comes from  $1 \ge {\frac{{\max\{a(p-1), 1\}+1}}{p+1}}$ and $1 \ge \frac{1}{p}$ for all $p \ge 4$ and $a \in [0, \frac{1}{2}]$  and  $\min \big\{ \epsilon ,  \epsilon^{\frac{1}{p}} ,  \epsilon^{\frac{{\max\{a(p-1), 1\}+1}}{p+1}} \big\}=\epsilon$  for $0<\epsilon<1$. 
\end{itemize}

\noindent
\textbf{Case 3 (locally nearly strongly convex):}   Using \eqref{nearlyconvex epi} and $\sigma_k \in \sigma_{\max}$,  we deduce\\
$
      \|s_k\| > \frac{1}{3} \min \bigg\{ \epsilon^{1-a} , \big(\frac{\epsilon}{L}\big)^{\frac{1}{p}} ,  \big(\frac{\epsilon}{\sigma_{\max}}\big)^{\frac{1}{p'-1}} \bigg\}. 
$
Let $\hat{\kappa}_4 :=\frac{1}{3} \min \bigg\{1,  L^{\frac{1}{p}},   \hat{C}_R^{-\frac{1}{p}}, \tilde{C}_R^{-\frac{1}{p+1}}   \bigg\} $. 
\begin{itemize}
    \item 
If $p$ is odd, $\sigma_{\max} =\hat{C}_R  \epsilon^{a(2-p)}$, then 
\begin{eqnarray}  
      \|s_k\| > 
     \underbrace{\frac{1}{3} \min \bigg\{ 1  , L^{\frac{1}{p}} , \hat{C}_R^{-\frac{1}{p}}\bigg\} }_{\ge \hat{\kappa}_4}
     \min \bigg\{ \epsilon^{1-a},   \epsilon^{\frac{1}{p}} ,\epsilon^{\frac{{a(p-2)+1}}{p}} \bigg\} \ge \hat{\kappa}_4 \epsilon^{1-a}      
    \label{nearly strongly convex 3 compare odd}
\end{eqnarray}
where the last inequality comes from  $1-a \ge \frac{{a(p-2)+1}}{p}$ and $1-a \ge   \frac{1}{p}$ for all $p\ge 3$  and $   \min \big\{  \epsilon^{1-a} , \epsilon^{\frac{1}{p}} ,\epsilon^{\frac{{a(p-2)+1}}{p}}  \big\}  = \epsilon^{1-a}   $ for $0<\epsilon <1$. 

\item If $p$ is even, $\sigma_{\max} :=\tilde{C}_R \max\big\{ \epsilon^{-1} ,   \epsilon^{a(1-p)}\big\} $, then
\begin{eqnarray}  
      \|s_k\| > \underbrace{
      \frac{1}{3} \min \bigg\{ 1, L^{\frac{1}{p}} ,  \tilde{C}_R^{-\frac{1}{p+1}} \bigg\}}_{\ge \hat{\kappa}_4}
       \min \bigg\{  \epsilon^{1-a} , \epsilon^{\frac{1}{p}} ,   \epsilon^{\frac{{\max\{a(p-1), 1\}+1}}{p+1}} \bigg\} \ge \hat{\kappa}_4 \epsilon^{1-a}    
    \label{nearly strongly convex 3 compare even}
\end{eqnarray}
where the last inequality comes from $1-a \ge {\frac{{\max\{a(p-1), 1\}+1}}{p+1}}$ and $ 1-a \ge \frac{1}{p}$ for all $p \ge 4$ and $a \in [0, \frac{1}{2}]$ and $   \min \big\{  \epsilon^{1-a} , \epsilon^{\frac{1}{p}} ,   \epsilon^{\frac{{\max\{a(p-1), 1\}+1}}{p+1}} \big\}  = \epsilon^{1-a}   $ for $0<\epsilon <1$.
\end{itemize}
\end{proof}

\begin{remark} \textbf{(The choice of bound for $a$)}
For odd $p\ge 3$, $1-a \ge \frac{{a(p-2)+1}}{p}$ is only true for $0\le a \le \frac{1}{2}$. For even $p\ge 4$,  $1-a \ge {\frac{{\max\{a(p-1), 1\}+1}}{p+1}}$ is only true for $0\le a \le \frac{1}{2}$. This rationale justifies the selection of $a \in [0, \frac{1}{2}]$. Further discussion regarding $a \geq a_b$ is provided in Appendix \ref{appendix a>1/2}.
\label{remark a>1/2}
\end{remark}

\begin{theorem}\textbf{(Function value decrease in successful iterations)}
\label{thm: Function Value Decrease in successful iter}
Let $f\in \C^{p, 1}(\R^n)$ with $p \ge 3$. Under Assumptions \ref{assumption Liptz}--\ref{assumption bounded hessian},
the mechanism of Algorithm \ref{arp + SoS algo} guarantees that, in any successful iteration $k \ge 0$ and for any $0< \epsilon <1$, 
\begin{eqnarray}
\text{Case 1 (odd $p$): } &\lambda_{\min}[H_k]\ge \epsilon^{a} , \qquad  & f(x_k) - f(x_k+s_k) >
\kappa_1\eta\epsilon^{\frac{p+1+a(p-2)}{p}},
      \label{Function Value Decrease convex odd}
 \\ \text{Case 1  (even $p$): } &\lambda_{\min}[H_k]\ge \epsilon^{a} , \qquad  &  f(x_k) - f(x_k+s_k) > 
\kappa_2\eta \epsilon^{\frac{p+2+\max\{1, a(p-1)\}}{p+1}},
  \label{Function Value Decrease convex even}
\\\text{Case 2: } & \lambda_{\min}[H_k] <0  ,   \qquad & f(x_k) - f(x_k+s_k)>  \frac{\eta \hat{\kappa}_3^2}{2 (\Lambda_2+1)}  \epsilon^2 ,
      \label{Function Value Decrease nonconvex}
\\\text{Case 3: }  &0 \le \lambda_{\min}[H_k]\le \epsilon^{a},  \qquad & f(x_k) - f(x_k+s_k) > 
\frac{\eta\hat{\kappa}_3^2}{2}\epsilon^{2 -a},
      \label{Function Value Decrease nearly strongly convex}
\end{eqnarray}
 where $\kappa_1, \eta, \hat{\kappa}_4^2, \hat{\kappa}_4^2$ and $\Lambda_2$ are iteration-independent and $\epsilon$-independent constants. 
\end{theorem}

\begin{proof}  Denote $E_k := f(x_k) - f(x_k+s_k)$.

\noindent
\textbf{Case 1 (locally strongly convex)}: Using \eqref{convex  SoS Taylor model}, we have
$$
\eta^{-1} E_k \ge f(x_k) - T_p(x_k,s_k) \underset{\eqref{temp mT 1}}{\ge} 
 \frac{\sigma_k }{p'}\|s_k\|^{p'} .
$$
According to Lemma \ref{Lemma Lower Bound for Step Size}, we have $  \|s_k\| > \frac{1}{2}\min\big\{\big(\frac{\epsilon}{L}\big)^{\frac{1}{p}}, \big(\frac{\epsilon}{\sigma_k}\big)^{\frac{1}{p'-1}} \big\} $. If $ \sigma_k\le \epsilon^{1-\frac{p'-1}{p}}
L^{\frac{p'-1}{p}}$ (i.e., for odd $p$, this is  $ \sigma_k\le L$), then we have $\|s_k\| >\frac{1}{2} \big(\frac{\epsilon}{L}\big)^{\frac{1}{p}}$. Consequently,
\begin{eqnarray}
\eta^{-1} E_k > \frac{\sigma_k}{2^{p'}p'} \bigg(\frac{\epsilon}{L}\bigg)^{\frac{p'}{p}} \ge \frac{\sigma_{\min}}{ 2^{p'}p'} \bigg(\frac{\epsilon}{L}\bigg)^{\frac{p'}{p}}
  \label{decrease both}
\end{eqnarray}
where the last inequality uses the lower bound $\sigma_k \ge \sigma_{\min}$ as defined in Algorithm \ref{arp + SoS algo}. 
Otherwise, if $ \sigma_k\ge \epsilon^{1-\frac{p'-1}{p}}L^{\frac{p'-1}{p}}$ (i.e., for odd $p$, this is  $ \sigma_k\ge L$), we have $\|s_k\| > \frac{1}{2}\big(\frac{\epsilon}{\sigma_k}\big)^{\frac{1}{p'-1}}$. Consequently, 
\begin{eqnarray}
 \eta^{-1} E_k> \frac{\sigma_k^{-\frac{1}{p'-1}}}{ 2^{p'}p'} \epsilon^{\frac{p'}{p'-1}}  \ge  \frac{1}{2^{p'}p'} {\sigma_{\max}}^{-\frac{1}{p'-1}} \epsilon^{\frac{p'}{p'-1}} 
 \label{dominating term}
\end{eqnarray}
\begin{itemize}
    \item 
If $p$ is odd, $p' = p+1$, $\sigma_{\max} = \hat{C}_R   \epsilon^{a(2-p)}$, then 
\begin{eqnarray}
 \eta^{-1} E_k > \frac{1}{2^{p'}p'} {\hat{C}_R}^{-\frac{1}{p}} \epsilon^{\frac{p+1+a(p-2)}{p}}.  
 \label{decrease odd}
\end{eqnarray}
Combining \eqref{decrease both} and \eqref{decrease odd}, we arrive that
$
\eta^{-1} E_k \ge  \frac{1}{2^{p'}p'} \min \big\{ \sigma_{\min} {L}^{\frac{p+1}{p}} {\epsilon}^{\frac{p+1}{p}},  \hat{C}_R^{-\frac{1}{p}}\epsilon^{\frac{p+1+a(p-2)}{p}} \big\}. 
$
Setting $\kappa_1 =   \frac{1}{2^{p'}p'} 
 \min \big\{ \sigma_{\min} {L}^{\frac{p+2}{p}}, \hat{C}_R^{-\frac{1}{p}}  \big\}$ gives in \eqref{Function Value Decrease convex odd}. 
\item If $p$ is even, $p' = p+2$, $\sigma_{\max} :=\tilde{C}_R \max\big\{ \epsilon^{-1} , \epsilon^{a(1-p)}\big\} $, then
\begin{eqnarray}
 \eta^{-1} E_k >  \frac{1}{2^{p'}p'} \tilde{C}_R^{-\frac{1}{p+1}} \epsilon^{\frac{p+2+ \max\{1, a(p-1)\}}{p+1}} .   
 \label{decrease even}
\end{eqnarray}
combining \eqref{decrease both} and \eqref{decrease even}, we arrive that
$
\eta^{-1} E_k \ge  \frac{1}{2^{p'}p'} \min \big\{ \sigma_{\min} {L}^{\frac{p+2}{p}} {\epsilon}^{\frac{p+2}{p}},  \tilde{C}_R^{-\frac{1}{p+1}}\epsilon^{\frac{p+2+\max\{1, a(p-1)\}}{p+1}} \big\}. 
$
Setting $\kappa_2 =   \frac{1}{2^{p'}p'} 
 \min \big\{ \sigma_{\min} {L}^{\frac{p+2}{p}}, \tilde{C}_R^{-\frac{1}{p+1}}  \big\}$ gives in \eqref{Function Value Decrease convex odd}\footnote{As indicated in \eqref{dominating term}, the dominating term for the function value reduction is $ f(x_k) - f(x_k+s_k) >  \frac{\eta}{2^{p'}p'} {\sigma_{\max}}^{-\frac{1}{p'-1}} \epsilon^{\frac{p'}{p'-1}}$. In the general AR$p$ framework \cite{cartis2020sharp, cartis2020concise}, we have $p'=p+1$ and $\sigma_{\max}$ is a constant independent of $\epsilon$. Therefore, $ f(x_k) - f(x_k+s_k) >  \mathcal{O}(\epsilon^{\frac{p+1}{p}})$ matching there result in \cite[Thm 3.7]{cartis2020concise}.
In the universal ARP$r$ framework \cite{cartis2019universal},  by setting $p'=p+2$, we have $\sigma_{\max} =  \mathcal{O}(\epsilon^{-\frac{1}{p}})$ \cite[Lemma 3.5]{cartis2019universal}. Therefore, $ f(x_k) - f(x_k+s_k) >  \mathcal{O}(\epsilon^{\frac{p+2}{p+1}}\epsilon^{\frac{1}{p(p+1)}})= \mathcal{O}(\epsilon^{\frac{p+1}{p}})$ matching there result in \cite[Thm 3.7]{cartis2019universal}.}. 
\end{itemize}

\noindent
\textbf{Case 2 (locally nonconvex):}  Using \eqref{nonconvex  SoS Taylor model}, we have 
\begin{eqnarray}
 \eta^{-1} E_k \ge f(x_k) - T_p(x_k,s_k) \underset{\eqref{temp mT 2}}{\ge} \frac{1}{2}\big(-\lambda_{\min}[H_k]+\delta\big)\|s_k\|^2
    \underset{\eqref{step low bound nonconvex}}{>}   \frac{\hat{\kappa}_3^2}{2}\big(-\lambda_{\min}[H_k]+\delta\big)^{-1} \epsilon^2 \ge \frac{\hat{\kappa}_3^2}{2}(\Lambda_2+1)^{-1} \epsilon^2
    \label{case 2 temp}
\end{eqnarray}
where the last inequality is obtained by Lemma \ref{lemma norm} and $0< \delta = \epsilon^a < 1$. 

\noindent
\textbf{Case 3 (locally nearly strongly convex):}  Using \eqref{nearly strongly convex  SoS Taylor model}, we have 
\begin{eqnarray}
 \eta^{-1} E_k \ge f(x_k) - T_p(x_k,s_k) \underset{\eqref{temp mT 3}}{\ge} \frac{\delta}{2}  \|s_k\|^{2}+ \frac{\sigma_k}{p'}\|s_k\|^{p'}   \underset{\eqref{step low bound nearly strongly convex}}{>}  \frac{\hat{\kappa}_4^2}{2} \epsilon^{a} \epsilon^{2-2a} = \frac{\hat{\kappa}_4^2}{2}  \epsilon^{2-a}. 
\label{case 3 temp}
\end{eqnarray}
\end{proof}

Based on Theorem \ref{thm: Function Value Decrease in successful iter}, we analyze how the function value reduction changes with $p$ and $\delta$. 
In Table \ref{table 1}, we give illustrations of the lower bound in \eqref{Function Value Decrease convex odd}--\eqref{Function Value Decrease nearly strongly convex} for different convexity cases.
We summarize the highlights as follows.

\begin{itemize}
    \item \textbf{For locally strongly convex iterations,} as $a$ decreases from $\frac{1}{2}$ to 0 (i.e., the lower bound of $\lambda_{\min}[H_k]$ increases from $\epsilon^{\frac{1}{2}}$ to $\mathcal{O}(1)$), if $p$ is odd, the function value reduction increases from $\mathcal{O}(\epsilon^{\frac{3}{2}})$ to $\mathcal{O}(\epsilon^{\frac{p+1}{p}})$. If $p$ is even, the function value reduction increases from $\mathcal{O}(\epsilon^{\frac{3}{2}})$ to $\mathcal{O}(\epsilon^{\frac{p+3}{p+1}})$. However, the region for Case 1 shrinks from $\lambda_{\min}[H_k] \ge \epsilon^{\frac{1}{2}}$ to $\lambda_{\min}[H_k] \ge c >0$ where $c >0$ is a constant with magnitude $\mathcal{O}(1)$. 
For odd $p$, when $a = 0$, the bound for the function value reduction in a locally strongly convex iteration is $\mathcal{O}(\epsilon^{\frac{p+1}{p}})$, matching the optimal complexity bound for a $p$th order method \cite{carmon2020lower, birgin2017worst, cartis2022evaluation, cartis2020sharp, cartis2020concise}. For even $p$, when $a = 0$, the bound for the function value reduction in a locally strongly convex iteration is $\mathcal{O}(\epsilon^{\frac{p+3}{p+1}})$. This arises from the choice to regularize our model with a $p'=p+2$ order of regularization instead of the $p'=p+1$ order proposed by \cite{carmon2020lower, birgin2017worst, cartis2022evaluation, cartis2020sharp, cartis2020concise}. The decision to use a $p'=p+2$ order of regularization when $p$ is even is motivated by the properties of the SOS Taylor model. As indicated by \cite{ahmadi2023higher}, even power in the regularization term is necessary to ensure that the model $m_p$ can be SOS-convex. 
    
    \item \textbf{For locally nearly strongly convex iterations,} conversely, as $a$ decreases from $\frac{1}{2}$ to 0 (i.e., the lower bound of $\lambda_{\min}[H_k]$ increases from $\epsilon^{\frac{1}{2}}$ to $\mathcal{O}(1)$), the function value reduction for Case 3 worsens from $\mathcal{O}(\epsilon^{\frac{3}{2}})$ to $\mathcal{O}(\epsilon^{2})$. The region for Case 3 expands from $0 \le \lambda_{\min}[H_k] \le \epsilon^{\frac{1}{2}}$ to $0 \le \lambda_{\min}[H_k] \le c$ where $c >0$ is a constant with magnitude $\mathcal{O}(1)$.
    \item \textbf{For locally nonconvex iterations,} regardless of how we choose $a \in [0, {\frac{1}{2}}]$, the function value reduction remains the same at $\mathcal{O}(\epsilon^2)$. This is expected, as we perturbed the second-order term by $-\lambda_{\min}[H_k]I_n + \delta I_n$. In locally nonconvex cases, the model $m_c$ starts to differ from $f$ (or $T_p$) from the second-order term of the Taylor expansion.  Therefore, due to the quadratic regularisation term added in the model in nonconvex iterations, this bound matches the accuracy order for high-order trust-region methods \cite{cartis2022evaluation}. 
    
    \item 
    In a locally strongly convex iteration with $a=0$, the magnitude of function value reduction increases for odd $p$ and even $p$, respectively. A higher-order method provides a larger function reduction on these iterations. However, note that for locally nonconvex iterations, the magnitude of the function value reduction remains unaffected by $p$ due to the term $\|s\|^2$.
\end{itemize}

\begin{table}[h!]
\caption{\small An illustration of the bound in \eqref{Function Value Decrease convex odd}--\eqref{Function Value Decrease nearly strongly convex} for different $a \in [0, \frac{1}{2}]$ in three convexity cases. }
\centering
 \begin{center}
\begin{tblr}{ |c||c|c|c| } 
 \hline[1.5pt]
\textbf{Odd $p \ge 3$} & \bb{$\lambda_{\min}[H_k]\ge \epsilon^{a}$} & \bb{$ 0< \lambda_{\min}[H_k]\le \epsilon^{a}$} &  \bb{$\lambda_{\min}[H_k] <0 $}  \\ 
 \hline[1pt]
 $a=0$ &  $\epsilon^{\frac{p+1}{p}}$ &  $\epsilon^2$  & $\epsilon^2$ \\ 
 $a=\frac{1}{3}$ &  $\epsilon^{\frac{4p+1}{3p}}$ &  $\epsilon^{\frac{5}{3}}$  & $\epsilon^2$ \\ 
 $a=\frac{1}{2}$ & $\epsilon^{\frac{3}{2}}$ & $\epsilon^{\frac{3}{2}}$  & $\epsilon^2$ \\ 
 \hline[1.5pt]
\end{tblr}
\end{center}

\centering
 \begin{center}
\begin{tblr}{ |c||c|c|c| } 
 \hline[1.5pt]
\textbf{Even $p \ge 4$}& \bb{$\lambda_{\min}[H_k]\ge \epsilon^{a}$} & \bb{$ 0< \lambda_{\min}[H_k]\le \epsilon^{a}$} &  \bb{$\lambda_{\min}[H_k] <0 $}  \\ 
 \hline[1pt]
 $a=0$ &  $\epsilon^{\frac{p+3}{p+1}}$ &  $\epsilon^2$  & $\epsilon^2$ \\ 
 $a=\frac{1}{3}$ &  $\epsilon^{\frac{4(p+1)+1}{3(p+1)}}$ &  $\epsilon^{\frac{5}{3}}$  & $\epsilon^2$ \\ 
 $a=\frac{1}{2}$ & $\epsilon^{\frac{3}{2}}$ & $\epsilon^{\frac{3}{2}}$  & $\epsilon^2$ \\ 
 \hline[1.5pt]
\end{tblr}
\end{center}
\label{table 1}
\end{table}

In the light of Theorem \ref{thm: Function Value Decrease in successful iter}, we propose practical choices for the parameter $a$ in Algorithm \ref{arp + SoS algo}.
\begin{remark} Let $0< \epsilon <1$. 
If $a = 0$ and $\delta = \mathcal{O}(1) = c$ where $c>0$ is constant, then
\begin{eqnarray*}
\text{locally strongly convex Iterations (odd $p$): }  &\lambda_{\min}[H_k]\ge c,&\qquad    
f(x_k) - f(x_k+s_k) > \mathcal{O}(\epsilon^{\frac{p+1}{p}}), 
\\\text{locally strongly convex Iterationss (even $p$): }  &\lambda_{\min}[H_k]\ge c,&\qquad    
f(x_k) - f(x_k+s_k) > \mathcal{O}(\epsilon^{\frac{p+3}{p+1}}), 
\\\text{Locally Nearly Convex Iterations: }& 0 \le \lambda_{\min}[H_k]\le c,& \qquad 
f(x_k) - f(x_k+s_k) > \mathcal{O}(\epsilon^{2}),
\\\text{Locally Nonconvex Iterations: }   &\lambda_{\min}[H_k]\le 0, & \qquad  f(x_k) - f(x_k+s_k) > \mathcal{O}(\epsilon^2).
\end{eqnarray*}

Let $a = {\frac{1}{2}}$  and $\delta = \epsilon^{\frac{1}{2}}$, then the order of complexity bound in \eqref{Function Value Decrease convex odd}, \eqref{Function Value Decrease convex even} and  \eqref{Function Value Decrease nearly strongly convex} equates,
\begin{eqnarray*}
\text{locally strongly convex Iterations: }  &\lambda_{\min}[H_k]\ge \epsilon^{\frac{1}{2}},&\qquad 
f(x_k) - f(x_k+s_k) > \mathcal{O}(\epsilon^{\frac{3}{2}}) ,
\\\text{Locally Nearly Convex Iterations: }& 0 \le \lambda_{\min}[H_k]\le \epsilon^{\frac{1}{2}},& \qquad 
f(x_k) - f(x_k+s_k) > \mathcal{O}(\epsilon^{\frac{3}{2}}),
\\\text{Locally Nonconvex Iterations: }   &\lambda_{\min}[H_k]\le 0, & \qquad  f(x_k) - f(x_k+s_k) > \mathcal{O}(\epsilon^2).
\end{eqnarray*}
\label{choice of delta}
\end{remark}

\subsubsection{Improved Function Value Decrease for $3 \le p \le 5$ in Locally Convex Iterations}
\label{sec Improved Function Value Decrease for p<5}

In this subsection, we demonstrate that for locally convex iterations with $3 \le p \le 5$, we can achieve an enhanced function value reduction per iteration. Since $m_c$ is SoS-convex, according to Definition \ref{def sos matrix} and Definition \ref{def sos convex}, we can deduce that $m_c$ is convex \cite{ahmadi2023higher, parrilo2000structured}. The idea of the proof is to use the convexity of $m_c$ (i.e., $\nabla_s^2 m_c(x_k, s) \succeq 0$  for all $s \in \R^n$) to derive a better bound for function value reduction in the locally convex iterations (i.e., Case 1). Note that the proof of the improved function value reduction bound in strongly convex iterations for $3 \le p \le 5$  is inspired by \cite{cartis2012evaluation}, although the analyses in \cite{cartis2012evaluation} are limited to the case $p=2$.

\begin{theorem}\textbf{(Function value decrease in successful iteration for $3 \le p \le 5$)}
\label{thm: p=2,3 Function Value Decrease in successful iteration}
Let $f\in \C^{p, 1}(\R^n)$ with $3 \le p \le 5$. Let $0< \epsilon <1$ and $\delta = \epsilon^a$ with $ a  \in [0, {\frac{1}{2}}]$.  Under Assumptions \ref{assumption Liptz}--\ref{assumption bounded hessian},
the mechanism of Algorithm \ref{arp + SoS algo} guarantees that, in successful iteration, for all $k \ge 0$,
\begin{eqnarray}
\text{Case 1 ($p=3,5$): } &\lambda_{\min}[H_k]\ge \epsilon^{a} , \qquad  &  f(x_k) - f(x_k+s_k) \ge\eta \hat{c} \epsilon^{a+\frac{2}{p}(1+a(p-2))},  \label{improved bound odd}
\\
\text{Case 1 ($p=4$): } &\lambda_{\min}[H_k]\ge \epsilon^{a} , \qquad  &  f(x_k) - f(x_k+s_k) \ge\eta \hat{c}_1  \epsilon^{a+\frac{2} {p+1}(1+ \max\{1, a(p-1)\})},
\label{improved bound even}
\\\text{Case 2: } & \lambda_{\min}[H_k] <0  ,   \qquad & f(x_k) - f(x_k+s_k) \ge 
 \frac{\eta \hat{\kappa}_3^2}{2 (\Lambda_2+1)}  \epsilon^2 ,  \notag
\\\text{Case 3: } &0 \le \lambda_{\min}[H_k]\le  \epsilon^{a},  \qquad & f(x_k) - f(x_k+s_k) \ge 
\frac{\eta\hat{\kappa}_4^2}{2}\epsilon^{2 -a}, \notag
\end{eqnarray}
as $\epsilon \rightarrow 0$, where $\hat{c}, \hat{c}_1, \eta$, and $\Lambda_2$ are iteration-independent and $\epsilon$-independent constants. 
\end{theorem}

\begin{proof} \textbf{(Proof for $p=3$)} The function value reductions for Cases 2 and 3 also follow similarly as in Theorem \ref{thm: Function Value Decrease in successful iter}. 
In this proof, we establish the function value reductions for Case 1 when $p=3$.  While the proof for $p=4,5$ follows a similar methodology as the proof for $p=3$, it is more intricate. Hence, we defer the proof for $p=4,5$ to Section \ref{sec proof for p=4,5}.

Firstly, we denote $f_k = f(x_k) \in \R$, $g_k = \nabla_x f(x_k) \in \R^n$, $ H_k  = \nabla_x^2 f(x_k) \in \R^{ n \times n}$ and $\T_k  = \nabla_x^3 f(x_k) \in \R^{ n^3}$. 
Then, $m_c$ and its derivatives have the following expression. 
\begin{eqnarray*}
 m_c(x_k, s) : &=& f_k + g_k^Ts + \frac{1}{2}H_k[s]^2 + \frac{1}{6}\T_k[s]^4+ \frac{\sigma_k}{4} \|s\|^4,
\\ \nabla_s m_c(x_k, s) : &=& g_k + H_k[s] + \frac{1}{2}\T_k[s]^2+ \sigma_k \|s\|^2s,
\\ \nabla_s^2 m_c(x_k, s) : &=& H_k +  \T_k[s]+ \sigma_k (\|s\|^2 I_n+ 2ss^T).
\end{eqnarray*} 
The function value reduction in each successful iteration is
\begin{eqnarray}
 \eta^{-1} E_k=  \eta^{-1}\bigg[  f(x_k) - f(x_k+s_k) \bigg] \ge f(x_k) - T_p(x_k,s_k) = -g_k^Ts_k -\frac{1}{2} H_k[s_k]^2 - \frac{1}{6}\T_k[s_k]^3.  
 \label{fun decrease temp}
\end{eqnarray}
Sine  $s_k = \argmin_{s \in \R^n} m_c(x_k, s)$, by first-order optimality condition $\nabla_s m_c(x_k, s)  =0$, we have
$
g_k =- H_k[s_k] - \frac{1}{2}\T_k[s_k]^2- \sigma_k \|s_k\|^2 s_k.
$
Note that an approximate first-order minimum of $m_c(x_k, s)$ is also allowed as discussed in Remark \ref{remark solver s_c}. Substituting the 1st order optimality condition into \eqref{fun decrease temp}, we have
\begin{eqnarray}
 \eta^{-1} E_k \ge \frac{1}{2} H_k[s_k]^2 + \frac{1}{3}\T_k[s_k]^3 + \sigma_k \|s_k\|^4. 
  \label{fun decrease temp 2}
\end{eqnarray}
Also, since $m_c(x_k, s)$ is SoS-convex with respect to $s$, we deduce that  $m_c(x_k, s)$ is convex, that is
\begin{eqnarray}
\bigg[H_k +  \T_k[s]+ \sigma_k (\|s\|^2 I_n+ 2 ss^T) \bigg] [u^2] \succeq 0, \qquad \forall  s, u \in \R^n. 
\label{convexity of mc}
\end{eqnarray}
By setting $s = \frac{3}{4}s_k$ and $u= \frac{2}{3} s_k$ in \eqref{convexity of mc}, we have $ \frac{4}{9} H_k [s_k]^2+  \frac{1}{3} \T_k[s_k]^3+  \frac{3\sigma_k}{4}    \|s_k\|^4  \ge 0$. Substituting this expression into \eqref{fun decrease temp 2}, and using $H_k \succeq \delta I_n$ for Case 1, we have
\begin{eqnarray*}
 \eta^{-1} E_k  \ge  \frac{1}{18} H_k[s_k]^2  +   \underbrace{   \frac{\sigma_k}{4}\|s_k\|^4}_{\ge 0}   \ge  \frac{\delta}{18} \|s_k\|^2 \underset{\eqref{step low bound convex}}{\ge}  \frac{\epsilon^a}{72} \min\bigg\{\bigg(\frac{\epsilon}{L}\bigg)^{\frac{2}{p}}, \bigg(\frac{\epsilon}{\sigma_{\max}}\bigg)^{\frac{2}{p}} \bigg\}
\end{eqnarray*}
where $\sigma_{\max} = \hat{C}_R  \epsilon^{a(2-p)}$. Set $\hat{c} =\frac{1}{72}  \min\{L^{-\frac{2}{p}}, \hat{C}_R^{-\frac{2}{p}}\}$, then
 $ 
 \eta^{-1} E_k  \ge \hat{c} \min \{\epsilon^{a+\frac{2}{p}}, \epsilon^{{a+\frac{2}{p}(1+a(p-2))}}\} = \hat{c} \epsilon^{{a+\frac{2}{p}(1+a(p-2))}}
$ which gives the desired bound. 
\end{proof}

We compare the improved rate of decrease in function value in Theorem \ref{thm: p=2,3 Function Value Decrease in successful iteration} to the rate given in Theorem \ref{thm: Function Value Decrease in successful iter}. 
In the locally convex case, the improved rate gives a large function value reduction for all choices of $a \in [0, {\frac{1}{2}}]$. 
In the light of Theorem \ref{thm: p=2,3 Function Value Decrease in successful iteration}, we propose practical choices for the parameter $a$ for $3 \le p \le 5$ in Algorithm \ref{arp + SoS algo}. 
\begin{remark} 
\label{special choice of delta}
Let $0< \epsilon <1$. Let $a = 0$ and $\delta = \mathcal{O}(1) = c$ where $c>0$ is constant, then
\begin{eqnarray*}
\text{locally strongly convex Iterations ($p=3,5$): }   &\lambda_{\min}[H_k]\ge c,&\qquad    
f(x_k) - f(x_k+s_k) > \mathcal{O}(\epsilon^{\frac{2}{p}}) ,
\\\text{locally strongly convex Iterations  ($p=4$): }   &\lambda_{\min}[H_k]\ge c,&\qquad    
f(x_k) - f(x_k+s_k) > \mathcal{O}(\epsilon^{\frac{4}{5}}),
\\\text{Locally Nearly Convex Iterations: }& 0 \le \lambda_{\min}[H_k]\le c,& \qquad 
f(x_k) - f(x_k+s_k) > \mathcal{O}(\epsilon^{2}), 
\\\text{Locally Nonconvex Iterations: }   &\lambda_{\min}[H_k]\le 0, & \qquad  f(x_k) - f(x_k+s_k) > \mathcal{O}(\epsilon^2).
\end{eqnarray*}
\end{remark}

\subsubsection{Proof for Improved Function Value Decrease for $p=4,5$ in Locally Convex Iterations}
\label{sec proof for p=4,5}

\begin{proof} \textbf{(Proof for $p=4,5$)}
We start by giving an expression for $m_c$ and its derivatives for $p=4,5$ in locally convex iterations (i.e., $\lambda_{\min}[H_k] \ge \delta $). Let $\T_k  = \nabla_x^3 f(x_k) \in \R^{ n^3}$, $\M_k : = \nabla_x^4 f(x_k) \in \R^{n^4}$ and $\N_k : = \nabla_x^5 f(x_k)  \in \R^{n^5}$. 
Then, $m_c$ and its derivatives have the following expression. Note that in the case of $p=4$, we have $\N_k =0$. 
\begin{eqnarray}
 m_c(x_k, s) : &=& f_k + g_k^Ts + \frac{1}{2}H_k[s]^2 + \frac{1}{3!}\T_k[s]^3 + \frac{1}{4!}\M_k[s]^4 + \frac{1}{5!}\N_k[s]^5+ \frac{\sigma_k}{6} \|s\|^6,
 \label{m5}
\\ \nabla_s m_c(x_k, s) : &=& g_k + H_k[s] + \frac{1}{2}\T_k[s]^2 + \frac{1}{6}\M_k[s]^3 + \frac{1}{24}\N_k[s]^4+ \sigma_k \|s\|^4s,
 \label{m5 grad}
\\ \nabla_s^2 m_c(x_k, s) : &=& H_k +  \T_k[s]+  \frac{1}{2}\M_k[s]^2 + \frac{1}{6}\N_k[s]^3+\sigma_k\|s\|^2 \bigg(\|s\|^2 I_n+ 4ss^T\bigg).
 \label{m5 hess}
\end{eqnarray} 
Let $ E_k:=    f(x_k) - f(x_k+s_k)$. In each successful iteration, $E_k$ satisfies
\begin{eqnarray}
 \eta^{-1}E_k  &\ge& f(x_k) - T_p(x_k,s_k) = -g_k^Ts_k -\frac{1}{2} H_k[s_k]^2 - \frac{1}{3!}\T_k[s_k]^3  - \frac{1}{4!}\M_k[s_k]^4 - \frac{1}{5!}\N_k[s_k]^5, \notag
\\ & =& \frac{1}{2} H_k[s_k]^2 + \frac{1}{3}\T_k[s_k]^3  + \frac{1}{8}\M_k[s_k]^4 + \frac{1}{30}\N_k[s_k]^5 + \sigma_k\|s_k\|^6
\label{m5 decrease}
\end{eqnarray}
where the second line  uses \eqref{m5 grad} and $\nabla_s m_c(x_k, s_k)^Ts_k =0$.

(\textbf{Scenario I}: assume $\N_k[s_k]^5 \ge 0$)  Since $m_c(x_k, s)$ is SoS-convex, we deduce that  $m_c(x_k, s)$ is convex.  Using $\nabla_s^2 m_c(x_k, s)[u]^2 \ge 0$ with $s = \frac{3}{4}s_k$ and $u= \frac{2}{3} s_k$, we have $ \frac{4}{9} H_k [s_k]^2+  \frac{1}{3} \T_k[s_k]^3+ \frac{1}{8}\M_k[s_k]^4 + \frac{1}{32}\N_k[s_k]^5 + \frac{45}{64}\sigma_k\|s_k\|^6\ge 0$. Substituting this expression into  \eqref{m5 decrease}, we have
\begin{eqnarray}
 \eta^{-1}E_k \ge   \bigg(\frac{1}{2}-\frac{4}{9} \bigg)\underbrace{H_k}_{\succeq \delta I_n}[s_k]^2  +\bigg( \frac{1}{30} -\frac{1}{32} \bigg) \underbrace{\N_k[s_k]^5}_{\ge 0} +  \bigg( 1 -\frac{45}{64} \bigg)\sigma_k \underbrace{\|s_k\|^6}_{\ge 0}\ge \frac{\delta }{18} \|s_k\|^2.
\end{eqnarray}
Note that the second inequality follows using the assumption of Scenario I ($N_k[s_k]^5 \ge 0$) and $H_k \succeq \delta I_n$. 

(\textbf{Scenario II}: assume $\N_k[s_k]^5 \le 0$) Using the convexity of $m_c$,  we have $    0 \le \nabla_s^2 m_c(x_k, s_k)[s_k]^2 = H_k[s_k]^2 +  \T_k[s_k]^3 + \frac{1}{2}\M_k[s_k]^4 + \frac{1}{6}\N_k[s_k]^5+ 5 \sigma_k\|s_k\|^6.$ For any $\tau >0$, this expression can be rearranged to give
\begin{eqnarray}
   \frac{\tau}{5}H_k[s_k]^2   \ge   \frac{\tau}{5} \bigg (- \T_k[s_k]^3 - \frac{1}{2}\M_k[s_k]^4- \frac{1}{6}\N_k[s_k]^5- 5 \sigma_k\|s_k\|^6\bigg).
   \label{temp tau}
\end{eqnarray} 
We can rewrite the second-order term in \eqref{m5 decrease} as $\big(\frac{1}{2}-\frac{\tau}{5}\big)  H_k[s_k]^2 + \frac{\tau}{5}H_k[s_k]^2 $ and substitute \eqref{temp tau} which gives
\begin{eqnarray}
 E_k  \ge \big(\frac{1}{2}-\frac{\tau}{5}\big) H_k[s_k]^2 + {\big(\frac{1}{3}-\frac{\tau}{5}\big)}\T_k[s_k]^3  + \big(\frac{1}{8}-\frac{\tau}{10}\big)\M_k[s_k]^4 + \big(\frac{1}{30}-\frac{\tau}{30}\big) \N_k[s_k]^5 + \sigma_k\big(1-\tau\big)\|s_k\|^6
\label{m5 decrease 2}
\end{eqnarray}
for all $\tau >0$.  Let $a, b \in \mathbb{R}$ and $s = a s_k $ and $u= b s_k$. Using convexity of $m_c$, we have 
\begin{eqnarray}
\nabla_s^2 m_c(x_k, s)[u]^2  = a^2 H_k [s_k]^2+  a^2b \T_k[s_k]^3+ \frac{a^2b^2}{2} \M_k[s_k]^4 +\frac{a^2b^3}{6}\N_k[s_k]^5 + 5a^2b^4\sigma_k\|s_k\|^6 \ge 0. 
\label{m5 hess coeff}
\end{eqnarray}
We can rewrite the second-order term in \eqref{m5 decrease 2} as $\big(\frac{1}{2}-\frac{\tau}{5} -a^2\big)  H_k[s_k]^2 + a^2 H_k[s_k]^2 $ and substitute \eqref{m5 hess coeff} which gives
\begin{eqnarray}
  \eta^{-1}E_k  \ge \big(\frac{1}{2}-\frac{\tau}{5} - a^2\big) H_k[s_k]^2 + {\big(\frac{1}{3}-\frac{\tau}{5} - a^2b\big)}\T_k[s_k]^3  + \big(\frac{1}{8}-\frac{\tau}{10}- \frac{a^2b^2}{2}\big)\M_k[s_k]^4 + 
 \notag
 \\ +\big(\frac{1}{30}-\frac{\tau}{30}-\frac{a^2b^3}{6}\big) \N_k[s_k]^5 + \sigma_k(1-\tau-5a^2b^4)\|s_k\|^6 
 \label{temp tau 2}
\end{eqnarray}
for all $\tau >0$, $a, b \in \R$. Set $\tau = \frac{4}{5}$, $a = \frac{26}{45} $ and $b = \frac{27}{52}$, we have $\frac{1}{3}-\frac{\tau}{5} - a^2b = 0$ and $\frac{1}{8}-\frac{\tau}{10}- \frac{a^2b^2}{2} = 0$ and \eqref{temp tau 2} becomes
\begin{eqnarray}
  \eta^{-1}E_k  \ge \frac{1}{162} \underbrace{H_k}_{\succeq \delta I_n}[s_k]^2 -\frac{7}{6240}\underbrace{\N_k[s_k]^5}_{\le 0}+ \frac{851}{10816}\sigma_k  \underbrace{\|s_k\|^6}_{\ge 0}\ge \frac{\delta}{162}\|s_k\|^2  \underset{\eqref{step low bound convex}}{\ge} \frac{\epsilon^a}{648}\min\bigg\{\bigg(\frac{\epsilon}{L}\bigg)^{\frac{2}{p}}, \bigg(\frac{\epsilon}{\sigma_{\max}}\bigg)^{\frac{2}{p'-1}} \bigg\}.
\end{eqnarray}
For $p=5$,  the rest of the analysis follows similarly as in the case of $p=3$.
For $p=4$, we have $\sigma_{\max} = \tilde{C}_R \max\big\{ \epsilon^{-1} ,   \epsilon^{a(1-p)}\big\}$. We set $\hat{c}_1 =  \frac{1}{648}\min \{ L^{-\frac{2}{p}}, \tilde{C}_R^{-\frac{2}{p+1}}\}$ and
$$
  \eta^{-1}E_k \ge \hat{c}_1 \min\big\{\epsilon^{a+\frac{2}{p}}, \epsilon^{a+\frac{2} {p+1}(1+ \max\{1, a(p-1)\})} \big\} = \hat{c}_1  \epsilon^{a+\frac{2} {p+1}(1+ \max\{1, a(p-1)\})}.
$$
\end{proof}


\subsection{Overall Complexity Bound}
\label{sec overall Complexity Bound}

We have discussed the complexity bounds for locally convex and nonconvex iterations. In this section, we provide the overall complexity bound for Algorithm \ref{arp + SoS algo}.

We define the sets of \textit{successful and unsuccessful iterations} for Algorithm \ref{arp + SoS algo} as $\mathcal{S}_m$ and $\mathcal{U}_m := [0:m]\backslash   \mathcal{S}_m$. Notice that  for $k \in \mathcal{S}_m$, $\rho_k \ge \eta$ and $x_{k+1} = x_k+ s_k$, while for  $ k \in  \mathcal{U}_m$, $\rho_k \le \eta$ and $x_{k+1} = x_k$. 
Lemma \ref{lemma Number of Iterations} describes the total number of iterations required to compute an $\epsilon$-approximate first-order minimizer using Algorithm \ref{arp + SoS algo}, as a function of successful iterations.

\begin{lemma} 
\label{lemma Number of Iterations}
\textbf{(Successful and unsuccessful adaptive-regularization iterations)}
    Suppose that Algorithm \ref{arp + SoS algo} is used and that $\sigma_k \le \sigma_{\max}$ for some $\sigma_{\max}>0$. Then, 
\begin{eqnarray}
     m \le |\mathcal{S}_m| \bigg(1+\frac{|\log \gamma_2|}{\log \gamma_1} \bigg)  + \frac{1}{\log \gamma_1} \log\bigg(\frac{\sigma_{\max}}{\sigma_0}\bigg).
     \label{no of iterations}
\end{eqnarray}
\end{lemma}

\begin{proof}
The update of regularization parameters in Algorithm \ref{arp + SoS algo}  gives that, for each $k$, 
$$
0\le \gamma_2 \sigma_k \le \max\{\bar{\sigma}_{k+1}, \sigma^{r}_{k+1}\} = \sigma_{k+1}, \quad k \in \mathcal{S}_m , \qquad \text{and} \qquad 
0\le \gamma_1 \sigma_j  \le \max\{\bar{\sigma}_{k+1}, \sigma^{r}_{k+1}\} = \sigma_{k+1}, \quad  k \in \mathcal{U}_m
$$
where $k \in [0:m]$. We inductively deduce that 
$  0 \le \gamma_2^{|\mathcal{S}_m |}\gamma_1^{|\mathcal{U}_m |}\sigma_0 \le \sigma_{k}.$ Therefore, using our assumption that $\sigma_k \le \sigma_{\max}$, we deduce that
$$
|\mathcal{S}_m | \log \gamma_2+|\mathcal{U}_m | \log \gamma_1 \le \log\bigg(\frac{\sigma_{\max}}{\sigma_0}\bigg). 
$$
The desired result follows from the inequality $m = |\mathcal{U}_m|+ |\mathcal{S}_m|$.
\end{proof}

Lemma \ref{lemma Number of Iterations} employs a similar technique to that in \cite[Lemma 2.3.1, Lemma 2.4.1]{cartis2022evaluation}. 
Note that Lemma \ref{lemma Number of Iterations} is independent of the form of the model $m_c$ and depends on the mechanism defined by the update of regularization parameters in Algorithm \ref{arp + SoS algo}.

\begin{lemma}
\label{lemma Bound on the number of successful iterations}
\textbf{(Bound on the number of successful iterations)} Let $f\in \C^{p, 1}(\R^n)$ with $p \ge 3$. Suppose that  Assumptions \ref{assumption Liptz}--\ref{assumption bounded hessian}  hold.  Let $f_{\text{low}}$ be a lower bound on $f(x)$ for $x \in \R^n$. Then, there exists a positive constant $\kappa_s$ such that Algorithm \ref{arp + SoS algo} requires at most
\begin{eqnarray}
  |\mathcal{S}_m|:= \kappa_s \frac{f(x_0)-f_{\text{low}}}{\epsilon^{2}}
  \label{bound for successful iterations}
\end{eqnarray}
successful iterations before an iterate $s_\epsilon$ is computed for which $\|g_k\| \le \epsilon$ for any $0< \epsilon <1$. 
\end{lemma}

\begin{proof}
According Theorem \ref{thm: Function Value Decrease in successful iter},  for all $k \in [0: (m-1)]$,
\begin{eqnarray}
f(x_k)-f(x_{k+1}) \ge \min \bigg\{ \kappa_1\eta\epsilon^{\frac{p+1+a(p-2)}{p}}, \quad \kappa_2\eta \epsilon^{\frac{p+2+\max\{1, a(p-1)\}}{p+1}}, \quad \frac{\eta \hat{\kappa_3}^2}{2 (\Lambda_2+1)}  \epsilon^2 ,  \quad \frac{\eta \hat{\kappa_4}^2}{2}\epsilon^{2 -a}\bigg\}.
\label{overall complexity}
\end{eqnarray}
Given that $p \ge 3$, $a \in [0, {\frac{1}{2}}]$, and $0 < \epsilon < 1$, by setting $\hat{\kappa}_s :=\min \bigg\{ \kappa_1\eta,  \kappa_2\eta, \frac{\eta \hat{\kappa_3}^2}{2 (\Lambda_2+1)},  \frac{\eta \hat{\kappa_4}^2}{2}\bigg\}$, we have $f(x_k)-f(x_{k+1}) \ge \hat{\kappa}_s\epsilon^2 $.
Suppose that $\|g_k\| > \epsilon$ for all $k \in [0: (m-1)]$. Then, 
\begin{eqnarray*}
f(x_0)-f_{\text{low}} \ge f(x_0)-f(x_k) = \sum_{k \in \mathcal{S}_m} \bigg[f(x_k)-f(x_{k+1})  \bigg] >  |\mathcal{S}_m|\hat{\kappa}_s\epsilon^2.
\end{eqnarray*}
By setting ${\kappa}_s :=\hat{\kappa}^{-1} $, we arrive at the bound \eqref{bound for successful iterations}. 

\end{proof}
Given Lemma \ref{lemma Bound on the number of successful iterations}, we are now ready to state the worst-case evaluation bound for Algorithm \ref{arp + SoS algo}. 

\begin{theorem} \textbf{(Complexity bound for Algorithm \ref{arp + SoS algo})} Let $f\in \C^{p, 1}(\R^n)$ with $p \ge 3$. Suppose that   Assumptions \ref{assumption Liptz}--\ref{assumption bounded hessian}  hold.  Then, there exists positive constants $\bar{\kappa}_s$,  $\kappa_c$,  $\kappa_d$ such that Algorithm \ref{arp + SoS algo} requires at most
$$
 \bar{\kappa}_s \frac{f(x_0)-f_{\text{low}}}{\epsilon^{2}}  + \kappa_c + \kappa_d 
 \log(\epsilon^{-1}) = \mathcal{O}(\epsilon^{-2}) 
$$
function evaluations (i.e., $f(x_k)$) and at most
$$
 \kappa_s \frac{f(x_0)-f_{\text{low}}}{\epsilon^{2}}  + 1
= \mathcal{O}(\epsilon^{-2})
$$
derivatives evaluations (i.e., $\nabla_x^jf(x_k)$ for $j \ge 1$) to compute an iterate $s_\epsilon$ such that $\|g_k\| \le \epsilon$ for any $0< \epsilon <1$ .
\label{thm: cqr complexity}
\end{theorem}

\begin{proof}
    The number of successful iterations needed to find an $\epsilon$-approximate first-order minimizer is bounded above by \eqref{bound for successful iterations}. The total number of iterations (including unsuccessful ones) can also be bounded above using Lemma \ref{lemma Number of Iterations}. 
    Using $\sigma_{\max} = \hat{C}_R \epsilon^{a(2-p)}$ for odd $p$ and    
    $\sigma_{\max} = \tilde{C}_R \max\big\{ \epsilon^{-1} ,   \epsilon^{a(1-p)}\big\} $ for even $p$ and   
   substituting \eqref{bound for successful iterations} into \eqref{no of iterations} yields
\begin{eqnarray*}
 m  \le \kappa_s \frac{f(x_0)-f_{\text{low}}}{\epsilon^{2}} \bigg(1+\frac{|\log \gamma_2|}{\log \gamma_1} \bigg) + \kappa_c  + \kappa_d \log(\epsilon^{-1}).
\end{eqnarray*}
Note that $\kappa_c :=  \frac{1}{\log \gamma_1} \log\big(\frac{\hat{C}_R}{\sigma_0}\big)$  and $\kappa_d  = \frac{a(p-2)}{\log \gamma_1} $ for odd $p$.  $\kappa_c :=  \frac{1}{\log \gamma_1} \log\big(\frac{\tilde{C}_R}{\sigma_0}\big)$ and $\kappa_d  = \frac{  \max \{ 1 , a(p-1) \}}{\log \gamma_1} $ for even $p$. 
Since Algorithm \ref{arp + SoS algo} uses at most one evaluation of $m_c$ per iteration and at most one evaluation of its derivatives per successful iteration (plus one evaluation of $m_c$ and its derivatives at the final iteration). 
By setting $ \bar{\kappa}_s = \kappa_s \big(1+\frac{|\log \gamma_2|}{\log \gamma_1} \big)  $, we deduce the desired conclusion. 
\end{proof}

\begin{remark}
The overall complexity bound is primarily influenced by the function value reduction in locally nonconvex cases. In these iterations, quadratic regularization is introduced to the nonconvex $m_c$ for SoS convexification. As a result, $m_c$ starts to deviate from $T_p$ or $f$ from the second-order term. The addition of the quadratic regularization term in nonconvex iterations gives the accuracy order that matches the accuracy order of the high-order trust-region methods \cite{cartis2022evaluation}.
\end{remark}

\subsubsection{Improved Complexity for Strongly Convex Functions}
\label{sec overall convex Complexity Bound}

For strongly convex functions $f(x)$, 
the evaluation complexity bound  for the $\epsilon$-approximate minimum value of the convex functions is $\mathcal{O}(\epsilon^{-\frac{1}{p}})$ for odd $p$, and the evaluation complexity bound of $\mathcal{O}(\epsilon^{-\frac{p+1}{p}})$ for even $p$. The global convergence rate for odd $p$ matches the bound by Nesterov in \cite[Thm 2]{Nesterov2021implementable}. 
It is worth noting that in \cite[Thm 5]{ahmadi2022complexity}, Ahmadi et al. derived the global convergence for odd orders $p \ge 3$ for the higher-order Newton method on a strongly convex objective function \cite[Thm 2]{Nesterov2021implementable}. In this section, we extend this result to include even orders $p \ge$ 2.  Also, note that neither \cite[Thm 5]{ahmadi2022complexity} nor \cite[Thm 2]{Nesterov2021implementable} discusses the adaptive regularization mechanism.  We provide the global convergence rate of  Algorithm \ref{arp + SoS algo} when minimizing a strongly convex objective function in Theorem \ref{rate of convergence convex}. To lay the groundwork, we introduce two technical lemmas.

\begin{lemma}(\cite[Lemma 3.2]{cartis2019universal})
 Let $D :=  \max_{1 \le j \le p} \big\{ \big(\frac{pp'}{j! \sigma_{\min}}\Lambda_j\big)^{\frac{1}{p'-j}}, 1\big\} >1$ be a constant independent of $k$, $\epsilon$ and $\delta$. Then, at each iteration of Algorithm \ref{arp + SoS algo}, we have
    \begin{eqnarray}
        \|s_k\|  \le D. 
        \label{assumption radius}
    \end{eqnarray}
\label{lemma radius}
\end{lemma}

\begin{lemma}
Let $f\in \C^{p, 1}(\R^n)$ with $p \ge 3$. Suppose that  Assumption \ref{assumption Liptz} holds and that $f$ is a strongly convex function such that $\nabla_x f(x) \succ \check{c} I_n \succ 0$ for some constant $  \check{c}>0$. The mechanism of Algorithm \ref{arp + SoS algo} guarantees that in a successful iteration, 
\begin{eqnarray}
        \label{technical convex result}
f(x_{k+1}) = f(x_k+s_k) \le  \min_{s \in \R^n } \bigg[f(x_k+s)  + \frac{L+D\sigma_{\max}}{p'}\|s\|^{p+1} \bigg].
\end{eqnarray}
where $D$ is defined in Lemma \ref{lemma radius}.  {For odd $p$, $\sigma_{\max} = \hat{C}_R \max \big\{1,\check{c}^{2-p} \}$, which is an iteration-independent and $\epsilon$-independent constant. For even $p$, $ \sigma_{\max} = \tilde{C}_R \epsilon^{-1}$ for any $0 < \epsilon <1$.} 
 \label{technical convex}
\end{lemma}

\begin{proof}  
This proof is adapted from \cite[Lemma 1]{Nesterov2021implementable} and uses similar techniques as those in \cite[Lemma 1]{Nesterov2021implementable}.
 {If $\nabla_x f(x) \succ \check{c} I_n \succ 0$, then in \eqref{general bound for sigma odd} and \eqref{general bound for sigma even}, we can replace $\delta$ with $ \check{c} $, which becomes an iteration-independent and $\epsilon$-independent constant. 
Consequently, for odd $p$, \eqref{odd bound a} becomes $  \sigma_k \le  \hat{C}_R \max \big\{1,\check{c}^{2-p} \}=\sigma_{\max}$.  For even $p$, \eqref{even bound a} becomes $  \sigma_k \le  \tilde{C}_R \max \big\{1, \epsilon^{-1}, \check{c}^{1-p} \} = \tilde{C}_R \epsilon^{-1}= \sigma_{\max}$ for $0 < \epsilon <1$. }

Assume $s_k$ is a successful iteration, then by \eqref{ratio test}, we have  $\eta^{-1}\big[  f(x_k) - f(x_k+s_k) \big] \ge f(x_k) - T_p(x_k,s_k)$  and $f(x_k+s_k)  \le f(x_k)$. Using these two inequalities, we deduce that 
 \begin{eqnarray*}
f(x_{k+1}) = f(x_k+s_k)&\le& ( 1- \eta^{-1})f(x_k) +  \eta^{-1}  f(x_k+s_k)  \le T_p(x_k,s_k)   
\\ &\le& T_p(x_k,s_k)+  \frac{\sigma_k}{p'}\|s_k\|^{p'}. 
 \end{eqnarray*}
By adding and subtracting $\frac{L}{p+1}\|s\|^{p+1}$, adding $\frac{\sigma_k}{p'}\|s\|^{p'}$ and \eqref{assumption radius}, we deduce
 \begin{eqnarray*}
 f(x_{k+1}) & \le & \min_{s \in \R^n } \bigg[\underbrace{T_p(x_k, s) - \frac{L}{p+1}\|s\|^{p+1}}_{\text{By \eqref{lip bound 1}, }\le f(x_k+s)} + \frac{L}{p+1}\|s\|^{p+1} + \frac{D\sigma_{\max}}{p'}\|s\|^{p+1} \bigg].
 \end{eqnarray*}
 \eqref{technical convex result} follows immediately from the fact that $\frac{1}{p'}< \frac{1}{p+1}$.
\end{proof}

\begin{theorem} 
Let $f\in \C^{p, 1}(\R^n)$ with $p \ge 3$. Suppose that  Assumptions \ref{assumption Liptz}--\ref{assumption bounded hessian}  hold.   Assume that $f$ is a strongly convex function such that $\nabla_x f(x) \succ \check{c} I_n \succ 0$ for some constant $  \check{c}>0$. Let the sequence $\{x_k\}_{k\ge 0}$ be generated by Algorithm \ref{arp + SoS algo}. Then, for $k\ge 0$, we have 
\begin{eqnarray}
        f(x_k)-f_*  \le  \nu  \bigg((p+1)^{\frac{1}{p}} + \frac{k-1}{p+1} \bigg)^{-p} \le \nu  \bigg(\frac{p+1}{k} \bigg)^{p}
        \label{linear rate}
\end{eqnarray}
where  $f_* =  \min_{x \in \R^n} f(x)$ and $\nu := (L+D \sigma_{\max})D^{p+1} >0$  is a constant independent of $k$ and $\epsilon$. 
 {For odd $p$, $\nu $  is an iteration-independent and $\epsilon$-independent constant. For even $p$, $\nu = (L+D \tilde{C}_R \epsilon^{-1})D^{p+1} =\mathcal{O}(\epsilon^{-1})$.}
\label{rate of convergence convex}
\end{theorem}

\begin{proof} This proof is adapted from \cite[Thm 2]{Nesterov2021implementable} and proof techniques follow closely from \cite[Thm 2]{Nesterov2021implementable}. 
Under the strongly convex assumption, $H_k \succeq \check{c} I_n \succ 0$ is always positive definite, the SoS Taylor model will take the form in \eqref{convex  SoS Taylor model}. 
We first prove \eqref{linear rate} for the case of $k=0$ and $x_k =0$. Using Lemma \ref{technical convex}, we derive
\begin{eqnarray}
 f(x_1)  \le \min_{s \in \R^n } \bigg[f(s) + \frac{(L+D \sigma_{\max})}{p+1}\|s\|^{p+1} \bigg] 
 \underset{\eqref{assumption radius}}{\le} f_* + \frac{(L+D \sigma_{\max})}{p+1}D^{p+1}. 
\label{first iter}
\end{eqnarray}
 {Let $\nu := (L+D \sigma_{\max})D^{p+1} >0$ with $\sigma_{\max}$ defined as in Lemma \ref{technical convex}.  Note that for odd $p$, $\nu $  is an iteration and $\epsilon$ independent constant and for even $p$, $\nu = (L+D \tilde{C}_R \epsilon^{-1})D^{p+1} =\mathcal{O}(\epsilon^{-1})$.}
From \eqref{first iter}, we have
\begin{eqnarray}
 f(x_1) - f_* \le  \frac{\nu}{p+1} < \nu  \bigg(\frac{p+1}{k} \bigg)^{p}
\label{first iter 1}
\end{eqnarray}
which satisfies \eqref{linear rate} for $k=1$ and $p \ge 3$. Further for $k \ge 1$, let  $x_* =  \argmin_{x \in \R^n} f(x)$, $f_* = f(x_*) = \min_{x \in \R^n} f(x)$ and $s_* = x_*-x_k$, using Lemma \ref{technical convex}, we have
 \begin{eqnarray*}
 f(x_{k+1}) &\le & \min_{s \in \R^n } \bigg[f(x_k+s) + \frac{(L+D \sigma_{\max})}{p+1} \| s\|^{p+1}  \bigg] 
\le  \min_{\alpha \in [0, 1] } \bigg[f\big[x_k+ \alpha s_* \big] + \frac{(L+D \sigma_{\max})}{p+1}\| \alpha s_*\|^{p+1}   \bigg] 
\\ &\underset{\eqref{assumption radius}}{\le} &   \min_{\alpha \in [0, 1] } \bigg[f\big[(1- \alpha)x_k+ \alpha (x_k+s_*) \big] + \frac{(L+D \sigma_{\max})}{p+1} D^{p+1} \alpha ^{p+1} \bigg]
 \\ &\le &   \min_{\alpha \in [0, 1] } \bigg[f(x_k) -  \alpha [f(x_k)-f_*] + \frac{\nu}{p+1} \alpha^{p+1} \bigg].
\end{eqnarray*}
The last inequality is obtained using the convexity of $f$.
The minimum of the above objective in $\alpha\ge 0$  is achieved for
$$
\alpha_* = \bigg( \frac{[f(x_k)-f_*]}{\nu } \bigg)^{\frac{1}{p}}  \le \bigg( \frac{[f(x_1)-f_*]}{\nu } \bigg)^{\frac{1}{p}}  \underset{\eqref{first iter 1}} {\le}\big(\frac{1}{p+1}\big)^{\frac{1}{p}} <1.
$$
Thus, we conclude that $  f(x_k+s_k) \le f(x_k) - \frac{p \alpha_*}{p+1} [f(x_k)-f_*]. $
Denoting $\Delta_k = f(x_k) -f_*$, we get the following estimate
$$
\Delta_k -\Delta_{k+1}= f(x_k)  - f(x_{k+1}) \ge \frac{p \alpha_*}{p+1} \Delta_k = C \Delta_k^{\frac{p+1}{p}}
$$
where $C = \frac{p}{p+1}  \nu ^{-1/p} $. 
Let $\mu_k  = C^p \Delta_k$. Following a similar recursive inequality analysis as in \cite[Thm 2]{Nesterov2021implementable}, we deduce that $$\mu_k^{-1} \ge \bigg( \mu_1^{-1/p} + (k-1)p^{-1} \bigg)^{p}
 \qquad \text{and}  \qquad
\mu_1^{-1/p} = C^{-1}   \Delta_1^{-1/p} \underset{\eqref{first iter 1}}{\ge}\frac{1} {p }(p+1)^{\frac{p+1}{p}}.
$$
Substituting these two inequalities into  $\mu_k  = C^p \Delta_k$, we have
\begin{eqnarray*}
\Delta_k &=& C^{-p} \mu_k \le  \bigg(\frac{p+1} {p }\bigg)^p \nu  \bigg( \mu_1^{-1/p} + (k-1)p^{-1} \bigg)^{-p} \le  \bigg(\frac{p+1} {p}\bigg)^p \nu  \bigg( \frac{1} {p }(p+1)^{\frac{p+1}{p}} + (k-1)p^{-1} \bigg)^{-p}
\\  &\le& \nu  \bigg((p+1)^{\frac{1}{p}} + \frac{k-1}{p+1} \bigg)^{-p} \ge \nu  \bigg(\frac{p+1}{k} \bigg)^{p}.
\end{eqnarray*}
\end{proof}

\begin{remark}
The assumptions for Theorem \ref{rate of convergence convex} closely resemble those in \cite[Theorem 2]{Nesterov2021implementable}. But Theorem \ref{rate of convergence convex} requires $f$ to be strongly convex, whereas \cite{Nesterov2021implementable} only requires $f$ to be convex. 
\end{remark}

\begin{remark} \label{Complexity bound for strongly convex objection functions}
\textbf{(Complexity bound for strongly convex objection functions.)}
The global convergence rates provided by Theorem \ref{rate of convergence convex} can be extended to establish a complexity bound for function and derivative evaluations when seeking the $\epsilon$-approximate minimum value of the convex functions.   {Similar to the approach in \cite[Thm 3.3]{grapiglia2020tensor} and \cite[Thm 3.4]{grapiglia2022tensor}, we bound \eqref{linear rate} by the tolerance $\epsilon$, 
$$ f(x_k)-f_*  \le \nu  \bigg(\frac{p+1}{k} \bigg)^{p} \le \epsilon$$
where $\nu$ is a constant for odd $p$ and $\nu$ has order $\mathcal{O}(\epsilon^{-1})$ for even $p$. Therefore, 
we obtained the evaluation complexity bound of $\mathcal{O}(\epsilon^{-\frac{1}{p}})$ for odd $p$, and the evaluation complexity bound of $\mathcal{O}(\epsilon^{-\frac{p+1}{p}})$ for even $p$.} While Theorem \ref{rate of convergence convex} specifies the convergence rate in successful iterations, by utilizing Lemma \ref{lemma Number of Iterations}, we can account for both successful and unsuccessful iterations to determine the total number of function and derivative evaluations required by Algorithm \ref{arp + SoS algo}.
\end{remark}

\section{Numerical Illustration of the Theoretical Bound}
\label{sec Numerics}

In this section, we explain how to implement the sum-of-squares-based paradigm \eqref{sos program} in Algorithm \ref{arp + SoS algo}. For implementation, we use the following lemma regarding SoS-convexity from \cite{ahmadi2013complete, kojima2003sums}.

\begin{lemma} \cite[Def 2]{ahmadi2023higher}
Let $y = [y_1, \dotsc, y_n]^T\in \R^n$ and $s = [s_1, \dotsc, s_n]^T\in \R^n$. Let $\R[s]^{n\times n}$ be the space of a $n \times n$ real polynomial matrix as defined in Definition \ref{def sos matrix}. We denote $\tilde{h}(s, y) := y^TH(s)y$. Note that $\tilde{h}(s, y)$ is a polynomial with $2n$ variables. Then, $H(s) \in \R[s]^{n\times n}$ is an SoS-matrix if and only if the polynomial $\tilde{h}(s, y)$  is a SoS polynomial\footnote{Note that the definition of  SoS polynomial is given in Definition \ref{def sos polynomial}. }. Namely, there exist polynomials $h_1, \dotsc, h_r: \R^{2n} \rightarrow \R$ such that $\tilde{h}(s, y) =\sum_{j=1}^r \Tilde{h}_j(s, y)^2$ for all $s,y \in \R^n$. 
\label{lemma sos convex}
\end{lemma}

Using Lemma \ref{lemma sos convex}, to determine whether $m_c$ in \eqref{sos program} is sos-convex, we need to verify whether $\tilde{h}(s, y) = y^T \nabla^2 m_c(x_k, s) y$ is sum-of-squares (SoS) polynomial.  
Established sum-of-squares parsers like YALMIP \cite{lfberg2004toolbox} and SOSTOOLS \cite{prajna2002introducing} provide functions such as \texttt{findsos} to determine whether a polynomial is SoS. Thus, to implement \eqref{sos program} in Algorithm \ref{arp + SoS algo}, we iteratively increase $\sigma_k$ until the polynomial $\tilde{h}(s, y) = y^T \nabla^2 m_c(x_k, s) y$  becomes SoS. Notably, $\tilde{h}(s, y)$ is a polynomial with $2n$ variables, $\{s_1, \dotsc, s_n, y_1, \dotsc, y_n\}$, and degree $p'$. Verifying whether a polynomial is SoS can be accomplished by solving a semidefinite program (SDP) of size polynomial in $n$ \cite{ahmadi2023higher}. Note that there are other methods and tools for sos-convex polynomial optimization problems, for instance, first-level Lasserre relaxation \cite{ahmadi2023higher, lasserre2001global}. Once $m_c$ is confirmed to be SoS convex, according to \cite[Corollary 2.5]{lasserre2008representation} and  \cite[Thm 3.3]{lasserre2009convexity}, we can determine its minimizer via another SDP of size polynomial in $n$.

\begin{figure}[!ht]
    \centering
    \includegraphics[width = 13cm]{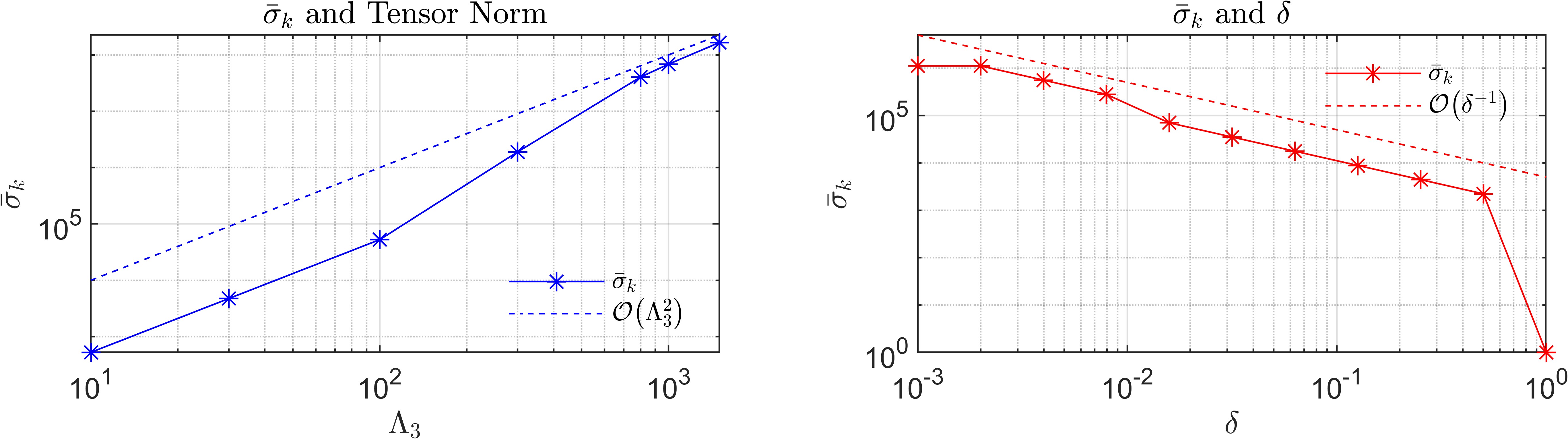} 
\caption{\small Numerical experiment conducted with $p=3$ and $n=2$. Here, $g_k \in \R^2$ and symmetric $H_k \in \R^{2 \times 2}$ contain randomly generated entries from $\mathcal{N}(0, 1)$. \textbf{Left}: Fix $ \delta = \mathcal{O}(1)$ and generate symmetric $\T_k \in \R^{2^3}$ with randomly generated entries, where the absolute maximum entries of $\T_k$ vary from $1$ to $10^3$. \textbf{Right}: Fix symmetric $\T_k \in \R^{2^3}$ containing randomly generated entries from $\mathcal{N}(0, 1)$ and vary $ \delta$ from $10^{-3}$ to $10^0$.}
    \label{fig: sigma and norm}
\end{figure}

We present preliminary numerical illustrations on how the size of the regularization parameter $\sigma_k$ changes with $\delta$ and the magnitude of $\nabla_x^j f(x_k)$.
\begin{itemize}
    \item \textbf{Investigating the relationship between $\bar{\sigma}_k$ and the magnitude of the tensor:} For a fixed $\delta$ and $p$, Corollary \ref{bound for sigma} indicates that the upper bound \eqref{sigma SoS bound} for $\bar{\sigma}_k$ scales like $(\sum_{j = 3}^p \frac{2 \Lambda_j}{R})^{p'-2}$. In the case of $p=3$, the upper bound simplifies to $\mathcal{O}\big({\Lambda_3}^2\big)$, where $\Lambda_3$ represents the largest absolute entry of the third-order tensor $\nabla_x^3 f(x_k) \in \R^{n^3}$. We plot the numerically obtained $\bar{\sigma}_k$ from \eqref{sos program} against the theoretical upper bound $\mathcal{O}({\Lambda_3}^2)$ in Figure \ref{fig: sigma and norm}. This plot illustrates that Corollary \ref{bound for sigma} provides a reasonable estimate for the size of $\bar{\sigma}_k$.
    \item 
\textbf{Investigating the relationship with $\bar{\sigma}_k$ and $\delta$:} According to Corollary \ref{bound for sigma}, the upper bound \eqref{sigma SoS bound} for $\bar{\sigma}_k$ scales like $\mathcal{O}(\delta^{3-p'})$. Specifically, for $p=3$, this yields $\mathcal{O}(\delta^{-1})$. Figure \ref{fig: sigma and norm} compares the theoretical bound $\mathcal{O}(\delta^{-1})$ with the actual value of $\bar{\sigma}_k$ obtained from \eqref{sos program}. The comparison indicates that our theoretical bound in \eqref{sigma SoS bound} provides a reasonable estimate for $\bar{\sigma}_k$.
\end{itemize}

\section{Conclusion}
\label{Conclusion}
In conclusion, this paper introduced an algorithmic framework that integrates the Sum of Squares (SoS) Taylor model with adaptive regularization techniques. 
In each iteration, we minimize a Sum of Squares (SoS) Taylor model which can be solved with polynomial cost per iteration. However, note that as $n$ increases, the size of SoS paradigm and its associated SDP can grow exponentially with the number of variables. This could cause challenges associated with the curse of dimensionality.
The function value reduction in each strongly convex iteration is $\mathcal{O}(\epsilon^{\frac{p+1}{p}})$ for odd $p$, and $\mathcal{O}(\epsilon^{\frac{p+3}{p+1}})$ for even $p$. Note that we have function value reduction for $p = 3, 5$ at $\mathcal{O}(\epsilon^{\frac{2}{p}})$ and for $p = 4$ at $\mathcal{O}(\epsilon^{\frac{4}{5}})$. For a general nonconvex function, the overall worst-case evaluation complexity bound is $\mathcal{O}(\epsilon^{-2})$.  {For strongly convex functions $f(x)$, 
the evaluation complexity bound  for the $\epsilon$-approximate minimum value of the convex functions is $\mathcal{O}(\epsilon^{-\frac{1}{p}})$ for odd $p$, and the evaluation complexity bound of $\mathcal{O}(\epsilon^{-\frac{p+1}{p}})$ for even $p$.}
Note that the bound does not address sub-problem costs, nor does it match the optimal complexity of the $p$th order optimization method in nonconvex iterations in \cite{carmon2020lower, cartis2022evaluation}. Yet, to the best of our knowledge, this is the first study that investigates both the global convergence and global complexity bounds for a tractable higher-order sub-problem ($p \ge 3$) for nonconvex smooth optimization problems.
This research also addresses a previously unanswered question in \cite{ahmadi2023higher} regarding the appropriate perturbation size for locally nonconvex iterations. 
We conducted a theoretical analysis of the perturbation parameter $\delta$ as a function of $\epsilon$, demonstrating its influence on the complexity bound of Algorithm \ref{arp + SoS algo} and giving some numerical guidance for its values. 

\smallbreak 
\noindent
\textbf{Acknowledgments}: We thank Karl Welzel for advice on this paper, particularly for his input on Remark \ref{remark bar H}.


\appendix
\section{Proof of Lemma \ref{lemma norm}}
\label{proof of lemma norm}
\begin{proof}
\textbf{Result 1)} Let $q(s) := \frac{1}{j!}\nabla_x^j f(x_k) [s]^j$ be the symmetric degree $j$ polynomial in $n$ variables, 
$$
q(s) = \frac{1}{j!} \nabla_x^j f(x_k) [s]^j = \frac{1}{j!} \sum_{i_1, \dotsc, i_j=1}^n   \big\{\nabla_x^j f(x_k)\big\}_{[i_1,\dotsc, i_j]} s_{[i_1]}\dotsc s_{[i_j]}
$$
where $1 \le i_1,\dotsc, i_j \le n$, $\big\{\nabla_x^j f(x_k)\big\}_{[i_1,\dotsc, i_j]}$ denotes the $\{i_1,\dotsc, i_j\}$th entry of the $j$th order tensor $ \nabla_x^j f(x_k)$ and $s_{[l]}$ denotes the $l$th entry of the vector $s \in \R^n$ for $1 \le l \le n$. Then the largest absolute value of the coefficients of $q(s)$ in a monomial basis can be expressed as 
$$
\|q(s)\|_{\infty^*} = \max_{1 \le i_1,\dotsc, i_j \le n}  \bigg\{\frac{1}{(j+1-\varsigma_j)!}\big\{\big|\nabla_x^j f(x_k)\big|\big\}_{[i_1,\dotsc, i_j]} \bigg\}\le  \max_{1 \le i_1,\dotsc, i_j \le n}  \big\{\big|\nabla_x^j f(x_k)\big|\big\}_{[i_1,\dotsc, i_j]}.
$$ 
where $1 \le \varsigma_j \le j$ denotes the number of distinct elements\footnote{Due to the symmetry of the tensor, we collect the coefficients corresponding to the same monomial basis element.} in $\{i_1,\dotsc, i_j\}$. 
We proved that $\|q(s)\|_{\infty^*}$ is smaller than the largest absolute value of the entry of $\nabla_x^j f(x_k)$.
Assume without loss of generality that the largest absolute value of the entry of $\nabla_x^j f(x_k)$ is attained at the $\{i_1^*,\dotsc, i_j^*\}$th entry, then 
$$
\|q(s)\|_{\infty^*} \le \big\{\big|\nabla_x^j f(x_k)\big|\big\}_{[i_1^*,\dotsc, i_j^*]} =\big|\nabla_x^j f(x_k)[e_{i_1^*}, \dotsc, e_{i_j^*}] \big| \le \big\|\nabla_x^j f(x_k) \big\|_{[j]} \underset{\text{Assumption } \ref{assumption bounded hessian}}{\le} \Lambda_j
$$ 
where $ e_{i_j^*} \in \R^n$ represents the unit vector with all entries $0$ except for the  $i_j^*$th entry, which is $1$. 

\noindent \textbf{Result 2)} For any $k \ge 0$, under Assumption \ref{assumption bounded hessian}, we have $\big|\lambda_{\min}[\nabla_x^2 f(x_k)]\big| \le \max_{\|v\|=1} |\nabla_x^2 f(x_k) [v]^2| \le \max_{\|v_1\|=\|v_2\|=1} |\nabla_x^2 f(x_k) [v_1, v_2] | =   \big\|\nabla_x^2 f(x_k) \big\|_{[2]} \le \Lambda_2$.
\end{proof}

\section{A Discussion for $a > {\frac{1}{2}}$}
\label{appendix a>1/2}

\textbf{For locally convex iterations,} the function value reduction in \eqref{Function Value Decrease convex odd} and \eqref{Function Value Decrease convex even} deteriorates as the value of $a$ increases. Opting for $a > \frac{1}{2}$ is unfavorable in this scenario. 

\noindent
\textbf{For the locally nonconvex case,} if $\frac{1}{2} < a < \frac{p'-2}{p'-3}$, then the same function value reduction bound in \eqref{Function Value Decrease nonconvex} remains. If $a> \frac{p'-2}{p'-3} > \frac{1}{2}$, then $\frac{{a(p-2)+1}}{p} \ge 1  \ge \frac{1}{p}$, $\max\{a(p-1), 1\} =a(p-1)$ and $\frac{{a(p-1)+1}}{p+1} \ge 1  \ge \frac{1}{p}$, therefore \eqref{step low bound nonconvex} becomes
$$
      \|s_k\| > 
      c_3 \epsilon^{\frac{a(p'-3)+1}{p'-1}} 
$$
where $c_3 = \frac{1}{3}\min\{ \hat{C}_R^{-\frac{1}{p}} , \tilde{C}_R^{-\frac{1}{p+1}}\}$. Consequently, \eqref{case 2 temp} becomes
$$
\eta^{-1} E_k > \frac{1}{2}\big(-\lambda_{\min}[H_k]+\delta\big)\|s_k\|^2 \ge \delta  \epsilon^{\frac{2a(p'-3)+2}{p'-1}} = \frac{c_3^2}{2} \epsilon^{a+\frac{2a(p'-3)+2}{p'-1}}. 
$$
Note that $a+\frac{2a(p'-3)+2}{p'-1} > 2+a > 2$ which is a worsen bound compared to \eqref{Function Value Decrease nonconvex}. 

\noindent
\textbf{For locally nearly strongly convex cases:} For odd $p \ge 3$ and $a> \frac{1}{2}$,  then $ 1-a \le \frac{{a(p-2)+1}}{p}$  and $ \frac{1}{p} \le  \frac{{a(p-2)+1}}{p}$. 
For even $p \ge 4$ and $a> \frac{1}{2}$,  then $\max\{a(p-1), 1\} = a(p-1)$, $1-a \le \frac{a(p-1)+1}{p+1}$ and $ \frac{1}{p} \le \frac{a(p-1)+1}{p+1}$. In both cases, the third term in \eqref{nearly strongly convex 3 compare even} and \eqref{Function Value Decrease convex even}  are dominant.  Consequently, $\|s_k\|> \frac{1}{3} {\hat{C}_R}^{-\frac{1}{p}}  \epsilon^{\frac{{a(p-2)+1}}{p}}$ for odd $p$ and  $\|s_k\|> \frac{1}{3} \tilde{C}_R^{-\frac{1}{p+1}} \epsilon^{\frac{a(p-1)+1}{p+1}}$ for even $p$.  Therefore \eqref{step low bound nearly strongly convex} becomes
$$
      \|s_k\| > 
      c_3 \epsilon^{\frac{a(p'-3)+1}{p'-1}} 
$$
where $c_3 = \frac{1}{3}\min\{ \hat{C}_R^{-\frac{1}{p}} , \tilde{C}_R^{-\frac{1}{p+1}}\}$. Consequently, \eqref{case 3 temp} becomes
$$
\eta^{-1} E_k > \delta  \frac{1}{2}\big(-\lambda_{\min}[H_k]+\delta\big)\|s_k\|^2 \ge \delta  \epsilon^{\frac{2a(p'-3)+2}{p'-1}} = \frac{c_3^2}{2} \epsilon^{a+\frac{2a(p'-3)+2}{p'-1}}. 
$$
Note that $\frac{2a(p'-3)+2}{p'-1}+a \ge \frac{3}{2}$ for $a >\frac{1}{2}$. The minimum of $\frac{2a(p'-3)+2}{p'-1}+a$ is attained at $a = \frac{1}{2}$ and decreases as $a$ increases. 

\scriptsize{
\bibliographystyle{plain}
\bibliography{sample.bib}
}

\end{document}